\documentclass[10pt]{article}  
\usepackage{amscd,amsfonts,amssymb,amscd,amsmath,latexsym}
\usepackage[hypertex,backref]{hyperref}
\makeatletter
\renewcommand{\subsubsection}{\@startsection
{subsubsection}
{1}
{0mm}
{0mm}
{0mm}
{\normalfont\normalsize\itshape}}
\makeatother
  
\newtheorem{theorem}{Theorem}[section] 
\newtheorem{propo}[theorem]{Proposition}
\newtheorem{lem}[theorem]{Lemma}
\newtheorem{ddd}[theorem]{Definition}

\newcommand{\forget}[1]{}

{\catcode`@=11\global\let\c@equation=\c@theorem}

\newcommand{\Z}{\mathbb{Z}}

\newcommand{\proof}{{\it Proof.$\:\:\:\:$}}

\newcommand{\R}{\mathbb{R}}

\newcommand{\Q}{\mathbb{Q}}

\newcommand{\C}{\mathbb{C}}

\newcommand{\Aut}{{\tt Aut}}

\newcommand{\cT}{\mathcal{T}}

\newcommand{\cH}{\mathcal{H}}

\newcommand{\cE}{\mathcal{E}}

\newcommand{\Hom}{{\tt Hom}}

\newcommand{\bc}{\mathbf{c}}

\newcommand{\coker}{{\tt coker}}
\newcommand{\id}{{\tt id}}

\newcommand{\nat}{\mathbb{N}}
\newcommand{\supp}{{\tt supp}}

\newcommand{\Ann}{{\tt Ann}}

\newcommand{\kpwcite}{\cite}
\def\imath{{i}}

\def\hB{\hspace*{\fill}$\Box$ \newline\noindent}

\def\hB{\hspace*{\fill}$\Box$ \\[0cm]\noindent}

 \newcommand{\cG}{\mathcal{G}}

\newcommand{\pr}{{\tt pr}}

\begin{document}
\title{$T$-duality for non-free circle actions\\
{\em \small Dedicated to Krzysztof Wojciechowski on his 50th Birthday}}
\author{Ulrich Bunke and Thomas Schick\footnote{Mathematisches Institut, Universit{\"a}t G{\"o}ttingen,
Bunsenstr. 3-5, 37073 G{\"o}ttingen, GERMANY,  bunke@uni-math.gwdg.de,
schick@uni-math.gwdg.de}}

\maketitle 
\abstract{
We study the topology of $T$-duality for pairs of $U(1)$-bundles and
three-dimensional integral cohomology classes over orbispaces.
}
\tableofcontents

\section{Introduction}

\subsection{From spaces to orbispaces}

\subsubsection{}

The concept of $T$-duality has its origin in string theory.
Very roughly speaking, it  relates one type of string theory on some target space 
with another type of string theory on a $T$-dual target space.  Some topological 
aspects of $T$-duality in the presence of $H$-fields were studied in
Bunke and Schick \kpwcite{bunkeschick045} (following earlier work  by
Bouwknegt, Mathai and Evslin \kpwcite{bem}, and others). In those preceeding investigations the main objects were pairs consisting of a $U(1)$-principal bundle 
and a three-dimensional integral cohomology class on its total space.
Here we could replace the notion of an $U(1)$-principal bundle by the equivalent notion of a free $U(1)$-space satisfying some slice condition. 

The main goal of the present paper is to extend the study of the topological aspects of $T$-duality to $U(1)$-spaces with finite stabilizers where we keep the slice condition. These spaces correspond to $U(1)$-bundles over orbispaces. 

\subsubsection{}

In order to deal properly with morphisms between orbispaces we will
use the more general language of topological stacks. Orbispaces are
particular topological stacks which admit an orbispace
atlas. Morphisms between orbispaces are required to be representable
maps. Our notion of an orbispace is a generalization of the notion of
a topological
space in the same spirit as the notion of an orbifold (see Moerdijk \kpwcite{mo}
for the definition of orbifolds which was motivating our definition of
orbispaces)
generalizes the notion of a smooth manifold.

 Topological $T$-duality is now
about pairs of $U(1)$-bundles in the category of orbispaces and 
three-dimensional cohomology classes in integral orbispace cohomology.
We will explain these notions at the appropriate places.

\subsubsection{}

Topological $T$-duality is the home for two different concepts. First
it is a relation on the set $P(B)$ of isomorphism classes of pairs
$(E,h)$ over a base space $B$, where $E\rightarrow B$ is a
$U(1)$-principal bundle and $h\in H^3(E,\Z)$ is an integral cohomology class
on the total space $E$ of the bundle. Secondly, $T$-duality denotes a
natural involution $T_B:P(B)\rightarrow P(B)$, which 
 associates to each pair a canonical isomorphism class of $T$-dual pairs.
In the present paper we generalize the definition of the $T$-duality relation as well as the
construction of canonical $T$-dual pairs (see \kpwcite{bunkeschick045}).
The main idea is to pass from orbispaces to spaces using a classifying space
functor. Once this functor is established the extension of the  results about the topology of $T$-duality of pairs from spaces to orbispaces is actually a formal matter.

\subsubsection{}

Another aspect of $T$-duality is the $T$-duality transformation
in twisted cohomology theories. It maps the twisted cohomology of the
total space of one $U(1)$-bundle to the twisted cohomology of its $T$-dual, where the twists are classified by the corresponding three-dimensional cohomology 
classes.
Of particular interest is the fact that under a $T$-admissibility assumption on the cohomology theory this transformation is an isomorphism.
In the present paper we discuss the generalization of this aspect to the orbispace case. In general it is a non-trivial matter to extend a cohomology theory to
the larger category of orbispaces. Of course, one could consider the Borel extension. In this case, where we again use the classifying space functor in order to pass from orbispaces to spaces, the generalization of the $T$-duality isomorphism is straight forward. On the other hand, having in mind the  example of $K$-theory, the Borel extension might not be the most interesting extension of the given generalized cohomology theory from topological spaces to orbispaces.

At the moment we do not know if the correct extension of twisted $K$-theory
to orbispaces is $T$-admissible.
 
\subsubsection{}

It is an amusing fact that the topology of $T$-duality
of $U(1)$-bundles over an orbispace as simple\footnote{Actually  the orbispaces $[*/\Gamma]$ are quite complicated. They are as complex as  the classifying space  $B\Gamma$.} as $[*/(\Z/n\Z)]$
(a point with the isotropy group $\Z/n\Z$) is already a non-trivial matter.
We will develop this example in detail.

This example serves as a building block of the more general example of
a Seifert bundle over a two-dimensional orbispace.
As an illustration we will calculate the $T$-dual of a Seifert bundle equipped with a three-dimensional cohomology class in terms of topological invariants.

\subsubsection{}

The problem of checking $T$-admissibility e.g. of twisted $K$-theory is
equivalent to the verification that the $T$-duality transformations
for all pairs over orbispaces of the form $[*/\Gamma]$ for all finite groups $\Gamma$ are isomorphisms. Currently we do not have explicit general results
about the topology of $T$-duality and the associated $T$-duality transformation
in this large class of examples.

\subsection{A detailed description of the contents}

\subsubsection{}

This paper is a continuation of \kpwcite{bunkeschick045}. In that paper we
introduced a contravariant set-valued homotopy invariant functor $P:spaces\rightarrow
sets$ which associates to each space $B$ the set of isomorphism
classes of pairs $(E,h)$ over $B$. Here $E\rightarrow B$ is a
$U(1)$-principal bundle and $h\in H^3(E,\Z)$. We have shown that the
functor can be represented by a space $R$ carrying a universal pair.
One of the main results was the determination of the homotopy type of
$R$. Consider the map $K(\Z,2)\times K(\Z,2)\rightarrow K(\Z,4)$ of Eilenberg-MacLane spaces given by the product of the canonical generators of the second cohomology of the two copies of $K(\Z,2)$. Then $R$ has  homotopy type of the homotopy fibre of this map.

\subsubsection{}

The notion of $T$-duality appeared first as a relation between
isomorphism classes of pairs. We then have shown that the universal
pair has a unique $T$-dual pair which determines and is determined by
its classifying map $T:R\rightarrow R$. This map induces a natural
transformation $T:P\rightarrow P$ which turns out to be two-periodic.

\subsubsection{}\label{not657}

The following short reformulation of the results of \kpwcite{bunkeschick045}
was suggested by the referee. It is close in spirit to the approach to $T$-duality for $U(1)^n$-principal bundles
via $T$-duality triples Bunke, Rumpf and Schick \kpwcite{bunkeschick045n}.
For two $U(1)$-principal bundles  $E\to B$ and $\hat E\to B$ let $E*\hat E\to B$ denote the fibrewise join.
It is a bundle with fibre $S^3$.
Let $\tilde P:spaces\to sets$ be the functor which associates to a space $B$ the set of isomorphism classes of triples $(E,\hat E,Th)$, where $Th\in H^3(E*\hat E,\Z)$ is a Thom class. 
Let $i:E\rightarrow E*\hat E$ be the natural inclusion map. Then $(E,\hat E,Th)\mapsto (E,i^*Th)$ defines a transformation $i:\tilde P\to P$. 
Using \kpwcite{bunkeschick045}, Thm. 2.16 one can show that this transformation is an isomorphism of functors. Under this isomorphism the $T$-duality transformation boils down to the involution
$T:\tilde P\rightarrow \tilde P$ given by $(E,\hat E,Th)\mapsto (\hat E,E,Th)$.
Note that this isomorphism $\tilde P\stackrel{\sim}{\to} P$ does not carry over to a corresponding result for $U(1)^n$-principal bundles if $n>1$, see  \kpwcite{bunkeschick045n}.

\subsubsection{}

There are various pictures of twisted cohomology theories.
In  \kpwcite{bunkeschick045} we descided to axiomatize those properties of twists and twisted cohomology theories which are used in connection with $T$-duality.

 In general, given a generalized cohomology theory
represented by some spectrum $E$ a twist of this cohomology theory over a space $B$ is something like a bundle of spectra with fibre $E$, or a presheaf of spectra with stalk $E$, depending on the framework.
The classification of twists is related to the classifying space $B\Aut(E)$ of the topological monoid of automorphisms of $E$. The twists considered in the present paper (as well as in the previous papers \kpwcite{bunkeschick045}, \kpwcite{bunkeschick045n}) are quite special and releated to the occurence of a map
$K(\Z,3)\rightarrow B\Aut(E)$ for cohomology theories like complex $K$-theory, $Spin^c$-cobordism theory,
or periodized real cohomology. In connection with $T$-duality the  restriction to this special sort of twists is crucial.

\subsubsection{}

In this setting,
twists should form a  functor $\cT:spaces\rightarrow groupoids$ such that the set of  isomorphism classes
of $\cT(B)$ is  in natural bijection with  $H^3(B,\Z)$, and such that the group of automorphisms of every $\cH\in \cT(B)$ is naturally isomorphic to $H^2(B,\Z)$.

In order to have an explicit model choose a realization of the Eilenberg-MacLane space $K(\Z,3)$. Then let $\cT(B)$ be the set of maps $B\to K(\Z,3)$. For two such maps $\cH,\cH^\prime$ let
$\Hom_{\cT(B)}(\cH,\cH^\prime)$ be the set of homotopy classes of homotopies from $\cH$ to $\cH^\prime$.

\subsubsection{}
In  \kpwcite{bunkeschick045} we have further introduced the notion of a $T$-admissible twisted
cohomology theory. It associates to a space $E$ and a twist $\cH\in \cT(E)$ the graded group $h(E,\cH)$. Twisted cohomology is functorial in both arguments.
If $u:\cH\rightarrow \cH^\prime$ is an isomorphism of twists, then
we have an induced map $u^*:h(E,\cH^\prime)\rightarrow h(E,\cH)$.
If $f:B^\prime\rightarrow B$ is a map of spaces, then
we have a functorial map $f^*:h(B,\cH)\rightarrow h(B^\prime,f^*\cH)$. It should furthermore admit an integration map for suitable oriented bundles. For details we refer to  \kpwcite{bunkeschick045}.

\subsubsection{}

Given a pair $(E,h)$ the class $h$
determines an isomorphism class $[\cH]$ of twists $\cH\in \cT(E)$. 
If $(\hat E,\hat h)$ is dual to $(E,h)$ and $[\hat \cH]=\hat h$, then the $T$-duality transformation
$$T:h(E,\cH)\rightarrow h(\hat E,\hat \cH)$$
is given by the following construction.
Note that there is a unique class $(E,\hat E,Th)\in \tilde P(B)$ such that $(E,h)\cong i(E,\hat E,Th)$ and
 $(\hat E,\hat h)\cong i\circ T(E,\hat E,Th)$ (see \ref{not657} for the notation).
Consider the fibre product
$$\begin{array}{ccccc}&&E\times_B\hat E&&\\
&p\swarrow&&\hat p\searrow\\
E&&&&\hat E\\
&\searrow&&\swarrow&\\
&&B&&\end{array}\ .$$
As explained in \kpwcite{bunkeschick045} the Thom class $Th$ determines an  isomorphism 
$u:\hat p^*\hat \cH\rightarrow p^*\cH$.
The $T$-duality transformation is defined as the composition
$$T:=\hat p_!\circ u^*\circ p^*\ .$$

\subsubsection{}

By definition, the twisted cohomology theory is $T$-admissible  
if the $T$-duality transformation is an isomorphism in the special
case where $B$ is a point.
In  \kpwcite{bunkeschick045} we have shown that $T$-admissibility implies, via a Mayer-Vietoris
argument, that the $T$-duality transformation is an isomorphism in
general.

\subsubsection{}

With these results our contribution consisted  in 
presenting an effective formalism and adding some precision and slight
generalizations
to the understanding of the topic as presented in \kpwcite{bem} or Mathai, Rosenberg \kpwcite{mr}.

In the present paper we develop a formalism which allows a
considerable generalization of $T$-duality. The spaces which were
suitable for $T$-duality in 
\kpwcite{bunkeschick045} were total spaces $E$ of principal
$U(1)$-fibrations $E\rightarrow B$. In particular, the spaces $E$ were free $U(1)$-spaces.

In the present paper we will relax this condition by admitting finite
stabilizers. In order to keep track of all information it turns out to
be necessary to consider the quotient $B:=[E/U(1)]$ as a topological
orbispace, i.e. as a proper topological stack on the category of topological
spaces which admit an orbispace atlas. For the language we refer to Heinloth \kpwcite{heinloth} and Noohi \kpwcite{math.AG/0503247},  but we will recall essential notions in Subsection \ref{lag}.
The brackets shall indicate that we consider the quotient as a stack and
not just as a space. The map $E\rightarrow [E/U(1)]$ is an atlas which 
represents $[E/U(1)]$ as a topological stack. Since $U(1)$ is compact, this stack is proper. The requirement that
$[E/U(1)]$ admits an orbispace atlas (note that $E\rightarrow [E/U(1)]$ is not an orbispace atlas) replaces the requirement of the existence of local trivializations in the case of principal bundles.

\subsubsection{}

Consider the simple example of the $U(1)$-stack $[U(1)/(\Z/n\Z)]$ (equipped with the trivial three-dimensional cohomology class)
which is actually  a space with a $U(1)$-action.
It will turn out that its canonical
$T$-dual is $U(1)\times [*/(\Z/n\Z)]$ (equipped with a non-trivial three-dimensional cohomology class). This stack is not equivalent to a
space. 
Therefore we are led to consider $U(1)$-bundles in the category of stacks
as the domain and the target of the canonical $T$-duality 
from the beginning.
By definition, a representable map $E\rightarrow B$ of topological
stacks is a $U(1)$-principal bundle, if it admits a fibrewise action
of $U(1)$, if in addition there is a $U(1)$-equivariant isomorphism  
$$\begin{array}{ccccc}E\times_B E&&\cong&& E\times U(1)\\
&\pr_1\searrow&&\pr_1 \swarrow\\
&&E&&\end{array} \ ,$$ where $U(1)$ acts on the second factors
(this means that $E\rightarrow B$ is a family of $U(1)$-torsors),
and if  for every map $T\rightarrow B$ with $T$ a space the induced map
$T\times_BE\rightarrow T$ has local sections.
 Note that $E\rightarrow [E/U(1)]$ is a $U(1)$-principal bundle in the category of stacks.



\subsubsection{}\label{not}
There are various equivalent ways to define the integral cohomology group 
$H^*(E,\Z)$ of a topological stack $E$.
One possibility is as the sheaf cohomology of the constant sheaf over
$E$ with fibre $\Z$. In the present paper we prefer to employ
classifying spaces. An atlas  $X\rightarrow E$ of the topological stack
gives rise to a topological groupoid $X\times_EX\Rightarrow X$ and
thus to a simplicial space $X^.$. Let $|X^.|$ denote its geometric
realization. If $E$ is an orbispace and $X$ is an orbispace atlas, then (see Proposition \ref{heq}) there is a natural isomorphism
$$H^*(E,\Z)\cong H^*(|X^.|,\Z)\ .$$

\subsubsection{}

A pair $(E,h)$ over a stack  $B$ will be a $U(1)$-principal bundle
$E\rightarrow B$ together with a class $h\in H^3(E,\Z)$.
Two pairs
$(E,h)$ and $(E^\prime,h^\prime)$ over $B$ are isomorphic if there exists
an isomorphism of $U(1)$-bundles $\phi:E\rightarrow E^\prime$ such
that $\phi^* h^\prime=h$. 

If $(E,h)$ is a pair over $B$, and $f:B^\prime\rightarrow B$ is a representable map of topological stacks, then we can define the pull-back
$f^*(E,h):=(f^*E,\tilde f^* h)$, where
$f^*E:=B^\prime\times_BE\rightarrow B^\prime$ is the induced $U(1)$-bundle,
and $\tilde f:f^*E\rightarrow E$ is the induced map.
This definition extends the functor
$P$ to a functor $P:(stacks,representable \:\:maps)\rightarrow sets$. Note that stacks form a two-category, and $P$ identifies two-isomorphic morphisms.

\subsubsection{}

Assume that $B$ is an orbispace, and let $Y\rightarrow B$ be an orbispace atlas of $B$. Let  $Y^.$ be the associated
simplicial space, and $|Y^.|$ be its geometric realization.
It turns out (Proposition \ref{heq})  that the homotopy type of
$|Y^.|$ is independent of the choice of $Y$ in a natural way.
In fact, if $i:Y^\prime\rightarrow Y$ is a refinement of orbispace atlases,
then $|i^.|:|(Y^\prime)^.|\rightarrow |Y^.|$
is a homotopy equivalence, where $i^.:(Y^\prime)^.\rightarrow Y^.$ is
the induced map of simplicial spaces. Furthermore, if $Y_1\rightarrow B$ is another orbispace atlas,
then the common refinement $Y\leftarrow Y\times_B Y_1\rightarrow Y_1$ is again an orbispace atlas.  

\subsubsection{}

A pair $(E,h)$ over $B$ gives rise to a pair
$(|X^.|,h)\in P(|Y^.|)$ as follows. Note that
$X:=Y\times_BE\rightarrow E$ is an orbispace atlas of $E$.
The natural map $X^.\rightarrow Y^.$ is a simplicial $U(1)$-bundle
which induces an ordinary $U(1)$-bundle $|X^.|\rightarrow |Y^.|$.
We can consider $h\in H^3(|X^.|,\Z)$.
Therefore given an orbispace atlas $Y\rightarrow B$ we obtain
a map
$$PA_Y:P(B)\rightarrow P(|Y^.|)\ .$$
The map is natural in $B$ and in the atlas $Y$ as follows.
Consider a representable map $f:B^\prime\rightarrow B$. Then we have the equality
$$PA_{Y^\prime}\circ f^*=|f^.|^*\circ PA_Y\ ,$$ where
$Y^\prime:=B^\prime\times_BY$ is the induced atlas of $B^\prime$, and
$f^.:(Y^\prime)^.\rightarrow Y^.$ is induced by
the natural map $Y^\prime\rightarrow Y$.

Consider now a refinement $i:Y^\prime\rightarrow Y$ of the orbispace atlas $Y\rightarrow B$. Then
we have the equality
$$|i^.|^*\circ PA_Y= PA_{Y^\prime}\ .$$

\subsubsection{}

The following theorem is the key to our generalization from spaces to orbispaces of the results
about $T$-duality of  pairs.

\begin{theorem}\label{main122}
If $B$ is an orbispace with orbispace atlas $Y\rightarrow B$, then
$PA_Y:P(B)\rightarrow P(|Y^.|)$
is an isomorphism.
\end{theorem}
This theorem will be proved in Section \ref{devo}.
The main intermediate result, Proposition \ref{ttz}, states that for a given orbispace atlas
$Y\rightarrow B$ the construction above on the level of $U(1)$-principal bundles provides an equivalence between the categories of $U(1)$-principal bundles
over $B$ and $|Y^.|$, where morphisms are homotopy classes of bundle isomorphisms.

\subsubsection{}

We use Theorem \ref{main122} and the naturality properties 
of the transformation
$PA_Y$ in order to extend the transformation
$T:P\rightarrow P$, which associates to an isomorphism class of pairs
a natural isomorphism class of $T$-dual pairs, from spaces to orbispaces.
Let $B$ be an orbispace and $Y\rightarrow B$ be an orbispace atlas.
\begin{ddd}\label{t321}
We define 
$T_B:P(B)\rightarrow P(B)$
by $$T_B:=PA_Y^{-1}\circ T_{|Y^.|}\circ PA_Y\ .$$
\end{ddd}
By Theorem \ref{main122} the map $T_B$ is well-defined.
It follows from the functorial properties of $PA_Y$ that
$T_B$ is independent of the choice of the orbispace atlas $Y\rightarrow B$.
It furthermore follows that the maps $T_B$ for all orbispaces
assemble to an automorphism of the functor $P$.

If $B$ is a space, then we can use the atlas $B\rightarrow B$. In this case
$T$ reduces to the original $T$ on spaces. Therefore our construction provides an extension of $T$ from spaces to orbispaces.
Since the original  $T$ on spaces is involutive, the same is true for
its extension to orbispaces.

\subsubsection{}

The second topic of the present paper is the $T$-duality
transformation in twisted cohomology.
To this end we first introduce the notion of a twisted cohomology
theory defined on orbispaces. Here we essentially repeat the axioms
formulated in \kpwcite{bunkeschick045} and add an axiom dealings with
two-isomorphisms. We show in Subsection \ref{tgw}  that every twisted
cohomology defined on spaces has a Borel extension to orbispaces. But
in general there might be different more interesting extensions
($K$-theory provides an example).

\subsubsection{}
Let us fix a twisted cohomology theory $h$ on orbispaces.
Given  two pairs $(E_i,h_i)$, $i=0,1$, which are $T$-dual (this is the $T$-duality relation, see \ref{rel}), we consider twists $\cH_i$ on $E_i$ classified by $h_i$. Then we define
a $T$-duality transformation $T:h(E_0,\cH_0)\rightarrow h(E_1,\cH_1)$ of degree one which is natural in $B$. We extend the notion of $T$-admissibility of a twisted cohomology theory to the orbispace case (Definition \ref{t1w}). If $h$ is $T$-admissible then
the $T$-duality transformation is an isomorphism (Theorem \ref{t2w}).

Compared with the case of spaces, in the case of orbispaces $T$-admissibility is much more complicated to check. The reason is that an orbispace can have a complicated local structure. At the moment we are not able to show that in the orbispace case twisted $K$-theory is $T$-admissible. But we shall see in Subsection \ref{tgw} that the Borel extension of a $T$-admissible twisted cohomology theory from spaces to orbispaces is again $T$-admissible.






 \subsubsection{}

The paper concludes with
the computation of the canonical $T$-duals in some instructive
examples in Section \ref{exex}.

\section{Some stack language}

\subsection{Topological stacks and orbispaces}\label{lag}

\subsubsection{}

In the present paper we consider stacks in topological spaces.
A stack is a sheaf of groupoids on this category. The sheaf conditions
are descend conditions for objects and morphisms with respect to open 
coverings of spaces. We refer to \kpwcite{heinloth}, \kpwcite{math.AG/0503247} for details.
Stacks form a two-category. 

The category of topological spaces 
is embedded into stacks by mapping a space $X$ to the sheaf
of sets $Y\mapsto \Hom(Y,X)$, and we consider a set as a groupoid with only identity morphisms. We  can and will consider spaces as stacks. 
This point of view is also reflected in our notation which uses the same type of letters for spaces and stacks.

\subsubsection{}\label{vvbb11}

We shall illustrate the stack notions in the example of quotient stacks.
Let $G$ be a topological group acting on a space $B$.
Then we can form the quotient stack $[B/G]$. It associates to a space $T$ the groupoid $[B/G](T)$ of pairs
$(P\to T,\phi)$, where $P\to T$ is a $G$-principal bundle and $\phi:P\rightarrow B$ is a $G$-equivariant map. The morphisms $(P\to T,\phi)\to (P^\prime\to T,\phi^\prime)$ are principal bundle isomorphisms
$P\rightarrow P^\prime$ which are compatible with the maps to $B$.
If $f:T^\prime\rightarrow T$ is a map of spaces, then $[B/G](f):[B/G](T)\rightarrow [B/G](T^\prime)$ is given by pull-back.

A $G$-equivariant map $h:B\rightarrow B^\prime$  induces a morphism of stacks
$h_*:[B/G]\rightarrow [B^\prime/G]$ by
$(P\to T,\phi)\mapsto (P\to T,h\circ \phi)$.

\subsubsection{}
A map $X\rightarrow Y$ between stacks
 is called representable if for each space $T$ and map $T\rightarrow Y$ the stack $T\times_YX$ is equivalent to a space.

\subsubsection{}\label{trewq}

Let us check that the map $h_*:[B/G]\rightarrow [B^\prime/G]$ in \ref{vvbb11} is representable. 
To this end we must calculate the fibre product $T\times_{[B^\prime/G]}[B/G]$ for a map
$f:T\rightarrow [B^\prime/G]$ and show that it is equivalent to a space. 
Let $f$ be given by $(P^\prime\to T,\phi^\prime)$.
We claim that $T\times_{[B^\prime/G]}[B/G]\cong(P^\prime\times_{B^\prime} B)/G$.
The map to $[B/G]$ is given by the pair
 $(P^\prime\times_{B^\prime} B\to  (P\times_{B^\prime} B)/G,\pr_2)$, and the map to $T$ is given by the composition $ (P^\prime\times_{B^\prime} B)/G\stackrel{\pr_1}{\to} P^\prime/G\cong  T$.

Let $S$ be a space. Then by definition of the fibre product of stacks an object in $(T\times_{[B^\prime/G]}[B/G])(S)$ is a triple
$(g,((P\to S),\phi),u)$, where $g:S\to T$ is an object of $T(S)$, i.e.   a map, $(P\to S,\phi)$ is an object of $[B/G](S)$, and $u:f(g)\rightarrow h(P\to S,\phi)$, i.e.
an isomorphism
$h:g^*P^\prime\rightarrow P$ of principal bundles such that $\phi^\prime\circ g^\sharp=\phi\circ h$,
where $g^\sharp:g^*P^\prime\rightarrow P^\prime$ is the induced map of total spaces.

The equivalence $(T\times_{[B^\prime/G]}[B/G])(S)\stackrel{\sim}{\to} ((P^\prime\times_{B^\prime} B)/G)(S)$ associates to $(g,((P\to S),\phi),u)$ the map
$S\rightarrow (P^\prime\times_{B^\prime}B)/G$ induced by
the $G$-equivariant map
$(g^\sharp\circ u^{-1},\phi):P\to P^\prime\times_{B^\prime}B$.

\subsubsection{}

A topological stack is a stack which admits an atlas.
An atlas of a stack $B$ is a representable map $X\rightarrow B$
from a space $X$ to $B$ which admits local sections. Here we say that
a map of stacks $X\rightarrow Y$ admits  local sections if for each
map
$T\rightarrow Y$ from a space $T$ to $Y$ each point  $y\in T$ has a
neighborhood $U\subset T$ such that there exists a map
$U\rightarrow X$ and a two-isomorphism from the composition 
$U\rightarrow X\rightarrow Y$ to the composition $U\rightarrow
T\rightarrow Y$.

A refinement of an atlas $X\rightarrow B$ is given by an atlas $X^\prime\rightarrow B$ and a diagram
$$\begin{array}{ccccc}
X^\prime&&\rightarrow&&X\\
&\searrow&&\swarrow&\\
&&B&&\end{array}\ .$$

\subsubsection{}

Let us check that the quotient stack $[B/G]$ considered in \ref{vvbb11} is topological.
We claim that $B\rightarrow [B/G]$ is an atlas. 

In order to see that this map is representable observe that
$B\cong [G/G]\times B\cong [(G\times B)/G]$, where in the last term $G$ acts on $G\times B$ by $h(g,b):=(gh^{-1},hb)$.
In order to see the first equivalence observe that $[G/G](S)$ is the groupoid of $G$-principal bundles with a section on $S$. This groupoid is connected and a set, hence equivalent to a one-point set. The second equivalence is induced by the $G$-equivariant map
$G\times B\rightarrow G\times B$, $(g,b)\mapsto (g,g^{-1}b)$, where the action of $G$ on the left $G\times B$ is given by $h(g,b):=(gh^{-1},b)$.
The map $B\cong [G\times B/G]\rightarrow [B/G]$ is now induced by the $G$-equivariant map
$\pr_2:G\times B \rightarrow B$. It is representable by \ref{trewq}.

Going through the definitions we see that the map $B\rightarrow [B/G]$ considered as an object of $[B/G](B)$ is given by $(G\times B\stackrel{\pr_2}{\to} B,\phi)$ with $\phi(g,b):=g^{-1}b$.

The existence of local sections can be seen as follows. Let $S\rightarrow [B/G]$ be a map given by a pair $(P\to S,\phi)$. Then we find a surjective map $f:A\to S$ such that $f^*P$ is trivial, i.e. admits an isomorphism $f^*P\cong G\times A$. The composition
$A\stackrel{a\mapsto (e,a)}{\longrightarrow} G\times A\cong f^*P\stackrel{f^\sharp}{\to} P\stackrel{\phi}{\to} B$ gives the required section.

\subsubsection{}

Given an atlas $X\rightarrow B$ we can define a topological groupoid
$$X\times_BX\Rightarrow X\ .$$
If $X^\prime\rightarrow X$ is a refinement, then we get an associated homomorphism of groupoids.

\subsubsection{}\label{p987}

In the case of the quotient stack $[B/G]$ with the atlas $B\rightarrow [B/G]$ this groupoid is the action groupoid $G\times B\Rightarrow B$, where the range and source maps are given by
$(g,b)\mapsto gb$ and $(g,b)\mapsto b$.

\subsubsection{}

A topological stack $B$ is called proper if the map of spaces
$$X\times_BX\rightarrow X\times X$$ is proper. This condition is independent of
the choice of the atlas.

\subsubsection{}\label{verrahgxdhjqw}

A topological groupoid $\cG^1\Rightarrow \cG^0$ is called {\'e}tale if the source and range maps $s,r:\cG^1\rightarrow \cG^0$ are {\'e}tale. 
A proper \'etale  topological groupiod $\cG^1\Rightarrow  \cG^0$ is called very proper if 
there exists a continuous function $\chi : \cG^0\to [0,1]$ such that
\begin{enumerate}
\item $r:\supp(s^*\chi)\to \cG^0$ is proper
\item  $\sum_{y\in r^{-1}(x)} \chi(s(y))=1$ for all $x\in \cG^0$.
\end{enumerate}
If $\cG$ is proper {\'e}tale such and  $\cG^0,\cG^1$ are locally compact spaces, then $\cG$  is automatically very proper. The existence of the corresponding cut-off functions has been shown e.g. in \cite[Prop. 6.11]{MR1671260}.

An orbispace atlas of a proper topological stack is an atlas $X\rightarrow B$ such that $X\times_BX\Rightarrow X$ is a very proper {\'e}tale topological groupoid.

We define a
topological  orbispace to be  a proper  topological stack which admits
an orbispace atlas. Our two-category of orbispaces $(orbispaces,
representable\:\: morphisms)$ has such orbispaces as objects and representable maps between orbispaces as one-morphisms.

 \subsubsection{}

We again consider quotient stack $[B/G]$ of \ref{vvbb11}. In view of \ref{p987} it is proper if and only if the action of $G$ on $B$ is proper, i.e. the map $G\times B\rightarrow B\times B$, $(g,b)\mapsto (gb,b)$, is proper. It is in addition  {\'e}tale if and only if $G$ acts with finite stabilizers. 

In particular, if $G$ is a discrete group acting properly on a locally compact space  $B$, then
$[B/G]$ is an orbispace.

\subsubsection{}

If $G$ is a finite group acting on the one-point space, then $[*/G]$ is an orbispace. If $G\rightarrow H$ is a homomorphism of finite groups, then we obtain
a map of stacks $[*/G]\rightarrow[*/H]$. It is a map of orbispaces
(i.e. representable) if and only if  the group homomorphism is injective.
In fact, in this case we can factor this map as
$[*/G]\cong [(G\backslash H)/H]\rightarrow [*/H]$, and the second map is prepresentable by \ref{trewq}.

\subsubsection{}

More generally, let $\cG:\cG^1\Rightarrow \cG^0$ be a topological groupoid acting on a space
$B$, i.e. there is a map $f:B\rightarrow \cG^0$ and an action
$B\times_{\cG^0}{} \cG^1\rightarrow B$ (the fibre product employs the range map $r:\cG^1\rightarrow \cG^0$).
Then we have the quotient stack
$[B/\cG]$. Its value on a space $X$ is given by the groupoid of pairs 
 $(P\rightarrow X,\phi)$ of locally trivial  
$\cG$-bundles $P\rightarrow X$ (see \kpwcite{heinloth}, Section. 3 for a definition) and maps $\phi:P\rightarrow B$ of $\cG$-spaces, and the  morphisms of the groupoid are the isomorphisms of such pairs.
There is a canonical map
$B\rightarrow [B/\cG]$ 
which is an atlas. Thus
$[B/\cG]$ is a topological stack.
If $\cG$ is proper and {\'e}tale and $B$ is locally compact, then $[B/\cG]$ is an orbispace.
In particular, we can apply this construction to the $\cG$-space $\cG^0$. We obtain the orbispace $[\cG^0/\cG]$ which is the classifying stack
for locally trivial $\cG$-bundles.

\subsection{Cohomology of orbispaces}

\subsubsection{}

Let $X\rightarrow B$ be an atlas of a topological stack and
$X\times_BX\Rightarrow X$ be the associated groupoid.
Then we obtain an associated simplicial space $X^.$ such that
$X^n:=\underbrace{X\times_B\dots\times_B X}_{n+1}$. By $|X^.|$ we denote its geometric realization.

A refinement $u:X^\prime\rightarrow X$ leads to a
map of simplicial spaces $u^.:(X^\prime)^.\rightarrow X^.$.
It further induces a map
$|u^.|:|(X^\prime)^.|\rightarrow |X^.|$ of realizations.

\subsubsection{}

In the present paper we heavily use the following fact
(which we learned from I. Moerdijk).
\begin{propo}\label{heq}
If $B$ is an orbispace, and $u:X^\prime\rightarrow X$ is a refinement of orbispace atlases of $B$, then
$|u^.|:|(X^\prime)^.|\rightarrow |X^.|$ is a weak homotopy equivalence of spaces.
\end{propo}
\proof
The category of sheaves (of sets) on the groupoid $X\times_BX\Rightarrow X$ is equivalent to the category of sheaves on $B$. In particular, the homomorphism of groupoids $$(X^\prime\times_BX^\prime\Rightarrow X^\prime)\rightarrow (X\times_BX\Rightarrow X)$$ induces an equivalence of categories of sheaves over groupoids.
In Moerdijk \kpwcite{moerdijk91} it is shown that
the category of sheaves on  $X\times_BX\Rightarrow X$ is equivalent
to the category of sheaves on the space $|X^.|$. 
If a map of spaces induces an equivalence of categories of sheaves, then it is a weak homotopy equivalence. This implies the result. \hB

\subsubsection{}

If $h(\dots)$ is some generalized cohomology theory then we can
extend this theory canonically to orbispaces. Given an orbispace $B$ we choose an orbispace atlas $X\rightarrow B$. Then we define
$$h(B):=h(|X^.|)\ .$$
This determines $h(B)$ up to natural isomorphisms (related to the various choices of the orbispace atlas).

If $f:B^\prime\rightarrow B$ is a representable map, then $X^\prime:=B^\prime\times_BX\rightarrow B^\prime$ is again an orbispace atlas. We obtain an induced morphism of groupoids
$(X^\prime\times_{B^\prime}X^\prime\Rightarrow X^\prime)\rightarrow
(X\times_BX\Rightarrow X)$, which induces a map of simplicial spaces $f^.:(X^\prime)^.\rightarrow X^.$, and eventually a map $|f^.|:|(X^\prime)^.|\rightarrow
|X^.|$ of geometric realizations.
The map
$f^*:h(B)\rightarrow h(B^\prime)$ is now given by
$|f^.|^*:h(|(X^\prime)^.|)\rightarrow h(|X^.|)$.

\subsubsection{}

Below we will apply this construction to integral cohomology
$h(\dots)=H(\dots,\Z)$.  
In order to distinguish the  construction described above from other
extensions of $h$ to orbispaces it will be called the Borel extension and denoted by $h_{Borel}$ (see also \ref{tgw}). This notation
is justified by its close relationship with the Borel extension of a cohomology theory to an equivariant cohomology theory.

\section{The $T$-duality relation}

\subsection{Thom classes and $T$-duality}\label{rel}

\subsubsection{}

Let $B$ be a topological stack.
We consider two $U(1)$-bundles
$E_i\rightarrow B$, $i=0,1$ over $B$ and let
$L_i\rightarrow B$ be the associated Hermitian vector bundles.
Let $S:=S(L_0\oplus L_1)\rightarrow B$ denote the unit-sphere bundle
in the sum of the two line bundles. Observe that the fibres of these bundles
are spaces since the corresponding projection maps to $B$ are representable.
We will denote points in the fibre of $S$ by
$(z_0,z_1)$, where $z_i\in L_i$ and $\|z_0\|^2+\|z_1\|^2=1$.
Then we have natural inclusions
$s_i:E_i\rightarrow S$ which identify
$E_i$ with the subsets $\{\|z_i\|=1\}$ for $i=0,1$, respectively.

\subsubsection{}

A Thom class for a three-sphere bundle $S\rightarrow B$ is a class
$Th\in H^3(S,\Z)$ which specializes to a Thom class of the three-sphere bundle
$|Y^.|\rightarrow |X^.|$ under the natural isomorphism $H^3(S,\Z)\cong H^3(|Y^.|,\Z)$ for some (and hence every) orbispace atlas $X\rightarrow B$, where $Y:=S\times_BX\rightarrow S$ is the induced atlas of $S$.

\subsubsection{}
Let $c_1(L_i)\in H^2(B,\Z)$ denote the first Chern classes of $L_i$.
As in the case of spaces  the three-sphere bundle $S\rightarrow B$ admits a Thom class if and only if
$c_1(L_0)\cup c_1(L_1)=0$ in $H^4(B,\Z)$.

\subsubsection{}
We now introduce the $T$-duality relation between pairs.
We consider classes $h_i\in H^3(E_i,\Z)$ for $i=0,1$ and the pairs
$(E_0,h_0)$ and $(E_1,h_1)$ over $B$.

\begin{ddd}\label{rel1}
We call the pairs
$(E_0,h_0)$ and $(E_1,h_1)$ $T$-dual if
there exists a Thom class $Th\in H^3(S,\Z)$ such that
$h_i=s_i^*Th$ for $i=0,1$, respectively.
\end{ddd}
This is the direct generalization of  \kpwcite{bunkeschick045},
Definition  2.9. 

\subsection{The $T$-duality transformation}

\subsubsection{}\label{rrg1}

In this subsection we assume that we have a twisted cohomology theory
defined on orbispaces. Thus given is a a functor of twists
$\cT:(orbispaces,representable \:\:maps)\rightarrow groupoids$ which satisfies the axioms listed
in \kpwcite{bunkeschick045}, Section 3.1 with spaces replaced by orbispaces.
As an additional datum we require that a two-isomorphism
$f\stackrel{\Phi}{\Rightarrow} f^\prime$ between maps $f,f^\prime:B^\prime\rightarrow
B$ induces an isomorphism of functors
$f^*\stackrel{\Phi_.}{\Rightarrow} (f^\prime)^*:\cT(B)\rightarrow \cT(B^\prime)$ in a
functorial way.

Furthermore, given is a bifunctor $h(\dots,\dots)$ which
associates to each pair $(B,\cH)$ of an orbispace $B$ and $\cH\in
\cT(B)$ a graded group $h(B,\cH)$, and which satisfies the axioms
listed again in \kpwcite{bunkeschick045}, Section 3.1.
 In addition we assume that
$f^*=\Phi_.^*\circ (f^\prime)^*:h(B,\cH)\rightarrow
h(B^\prime,f^*\cH)$ for two-isomorphic morphisms using the
notation above.

We require that the integration map $g_!:h(B^\prime,g^*\cH)\rightarrow h(B,\cH)$ is
defined for representable proper maps $g:B^\prime\rightarrow B$ which
are $h$-oriented. By definition, the datum of an $h$-orientation of $g$ is
equivalent to a compatible choice of $h$-orientations of the induced maps
of spaces $T\times_BB^\prime\rightarrow T$ for all maps $T\rightarrow B$, where
$T$ is a space.

\subsubsection{}
We consider an orbispace $B$.
Let $(E_0,h_0)$ and
$(E_1,h_1)$ be pairs over $B$ and $Th\in H^3(S,\Z)$ be a Thom class such that
$s_i^*Th=h_i$. We choose a twist $\cH\in \cT(S)$ such that $[\cH]=Th$. Then we define the twists $\cH_i:=s_i^*\cH\in \cT(E_i)$ for $i=0,1$.
In the present section we define the $T$-duality transformation
$$T_0:h(E_0,\cH_0)\rightarrow h(E_1,\cH_1)\ .$$

\subsubsection{}

We consider the two-torus bundle $F:=E_0\times_BE_1\rightarrow B$.
The map
$$ F\ni (z_0,z_1)\rightarrow (\frac{1}{\sqrt{2}}z_0,\frac{1}{\sqrt{2}}z_1) \in S$$ defines embedding which  gives rise to a decomposition
$$S\cong S_0\cup_F S_1\ ,$$
where $$S_i:=\{(z_0,z_1)\in S\:|\: \|z_i\|\ge \|z_{1-i}\|\}\ .$$

\subsubsection{}

The composition
$s_0\circ \pr_0 :F\rightarrow S$ is homotopic to the inclusion by the homotopy
$$(z_0,z_1)\mapsto  (\sqrt{1-\frac{t}{2}} z_0 ,\sqrt{\frac{t}{2}}z_1)\ , t\in [0,1]\ .$$
Similarly, $s_1\circ \pr_1$ is homotopic to the inclusion.
These homotopies give rise to isomorphism classes of isomorphisms of twists
$$v_i:\cH_{|F}\stackrel{\sim}{\rightarrow} \pr_i^*\cH_i\ .$$

\subsubsection{}

\begin{ddd}\label{tde2}
We define the $T$-duality transformations
$$T_i:h(E_i,\cH_i)\rightarrow h(E_{1-i},\cH_{1-i})$$
as the compositions
$$T_i:=(\pr_{1-i})_!\circ (v_{1-i}^{-1})^*\circ v_i^*\circ \pr_i^*\ .$$
\end{ddd}
Here it is essential to use the transformation
$(v_{1-i}^{-1})^*\circ v_i^*:\pr_{1-i}^*\cH_{1-i}\rightarrow \pr_i^*\cH_i$.
With other choices we can not expect that the maps $T_i$ become isomorphisms 
for $T$-admissible cohomology theories.










\subsection{$T$-admissible cohomology theories}

\subsubsection{}\label{admis1}

Let $\Gamma$ be a finite group, and choose two characters
$\chi_0,\chi_1:\Gamma\rightarrow U(1)$. We consider the stack
$B:=[*/\Gamma]$ and 
the bundles $E_i:=[U(1)/_{\chi_i}\Gamma]\rightarrow [*/\Gamma]$,
where $\Gamma$ acts on $U(1)$ by $\chi_i$ (this is indicated by the subscript), $i=0,1$.
We further consider classes $h_i\in H^3(E_i,\Z)$ such that
$(E_0,h_0)$ and $(E_1,h_1)$ are $T$-dual according to Definition \ref{rel1}. This is a non-trivial
condition
as we shall see later in \ref{gp1}.

\begin{ddd}\label{t1w}
Following
\kpwcite{bunkeschick045}, Definition 3.1,2 we call a twisted cohomology
theory $h(\dots,\dots)$ on orbispaces $T$-admissible if the $T$-duality
transformations $T_i$ are isomorphisms for all examples of the type
described above (i.e. for all choices finite groups $\Gamma$, pairs of characters $\chi_0,\chi_1$, and choices of the classes $h_i$).
\end{ddd}

\subsubsection{}

If the cohomology theory is $T$-admissible then the property that the $T$-duality transformation is an isomorphism can be extended to the large class
of base orbispaces $B$ which are build by glueing the local examples of the form
$[*/\Gamma]$. The argument is based on the Mayer-Vietoris sequence.

We call an orbispace $B$ finite if it has a finite filtration
$$\bigsqcup_i^{finite} [*/\Gamma_{i,0}]=B^0\subset
B^1\subset\dots\subset B^r=B$$ such that there exists
cartesian diagrams
\begin{equation}\label{eq14}\begin{array}{ccc}
S^{n_\alpha-1}\times[*/\Gamma_\alpha]&\rightarrow &B^{\alpha-1}\\
\downarrow&&\downarrow\\
D^{n_\alpha}\times
[*/\Gamma_\alpha]&\stackrel{i_\alpha}{\rightarrow}&B^\alpha
\end{array}\end{equation}
for $n_\alpha\in \nat$ and appropriate finite groups $\Gamma_\alpha$,  where
the $i_\alpha$ are representable and induce inclusions of open substacks
$(D^{n_\alpha}\setminus S^{n_\alpha-1})\times
[*/\Gamma_\alpha]\rightarrow B^\alpha$ (see \kpwcite{heinloth},
Definition 2.8), and $D^{n_\alpha}\times [*/\Gamma_\alpha]\sqcup B^{\alpha-1}\rightarrow B^\alpha$ is surjective.

For example, if  $M$ is a compact smooth manifold on which a compact group $G$ acts with finite stabilizers, then $[M/G]$ is a finite orbispace. In fact, $M$ admits a $G$-equivariant triangulation (by $G$-simplices of the form $\Delta^k\times G/H$ with $H\subset G$ a finite subgroup). Using this triangulation
we obtain the required filtration of $[M/G]$.
We expect that compact orbifolds in the sense of \kpwcite{mo} are finite orbispaces.

\subsubsection{}

\begin{theorem}\label{t2w}
Assume that the twisted cohomology theory is $T$-admissible.
Let $B$ be a finite orbispace, and let $(E_0,h_0)$ and $(E_1,h_1)$ be pairs
over $B$ which are $T$-dual to each other. Then the $T$-duality
transformations \ref{tde2} are isomorphisms.
\end{theorem}
\proof
This theorem is proved using induction over the number of cells of $B$
and the Mayer-Vietoris sequence in the same way as
\kpwcite{bunkeschick045}, Thm. 3.13. \hB

Using the method of the proof of Proposition \ref{rt51} we could
weaken the finiteness condition.


\subsubsection{}

It is natural to expect that an appropriate extension of twisted Atiyah-Segal $K$-theory to orbispaces is $T$-admissible. At the moment we do not
have a proof. In the following Subsection \ref{tgw} we provide examples of $T$-admissible cohomology theories.

\subsection{Borel-$K$-theory as an admissible cohomology theory on orbispaces}\label{tgw}

\subsubsection{}

The goal of the present subsection is to show that every twisted
cohomology theory defined on spaces and satisfying the list of axioms
stated in \kpwcite{bunkeschick045}, Section 3.1, admits an extension to
orbispaces by a Borel construction. For a demonstration we use $K$-theory.
We shall see that the Borel extension of a $T$-admissible twisted cohomology theory is again $T$-admissible.

\subsubsection{}

Note that in the case of $K$-theory the Borel construction is probably
not the most interesting extension to orbispaces. A better extension is
provided by the construction of Tu, Xu and Laurent \kpwcite{lht}.



\subsubsection{}

An extension of  a twisted cohomology theory from spaces to orbi\-spaces 
consists of an extension of the notion of a twist from spaces to orbi\-spaces, and then of the extension of the cohomology functor itself.

We start with the discussion  of twists. In this subsection we
will assume that we are given a functor $\cT$ on spaces which
associates to each space $B$ the groupoid of twists $\cT(B)$ (Note
that in general twists form a two-category. Here we adjust the notion by 
identifying isomorphic isomorphisms.)

\subsubsection{}

We now extend twists to orbispaces.
\begin{ddd}\label{tws}
A twist of an orbispace $B$ is given by an orbispace atlas $X\rightarrow B$ and a twist
$\cH\in \cT(|X^.|)$. A morphism of twists $\cH\rightarrow \cH^\prime$, where $\cH\in \cT(|X^.|)$ and $\cH^\prime\in\cT(|(X^\prime)^.|)$,  is given by
a common refinement $Y\rightarrow B$ of the orbispace atlases $X$ and $X^\prime$ and a morphism
$\phi:u^*\cH\rightarrow (u^\prime)^*\cH^\prime$,
where $u:|Y^.|\rightarrow |X^.|$ and
$u^\prime:|Y^.|\rightarrow |(X^\prime)^.|$ are the induced maps.
\end{ddd}
We identify morphisms which become equal on a common refinement of orbispace atlases.
In this way we associate to each orbispace $B$ a category of
twists $\cT(B)$.

\subsubsection{}\label{funcrf}

Let $f:B^\prime\rightarrow B$ be a morphism of orbispaces, i.e. a representable map of stacks.
Then we define the pull-back $f^*:\cT(B)\rightarrow \cT(B^\prime)$ as
follows.
If $X\rightarrow B$ is an orbispace atlas then we get an orbispace atlas
$X^\prime:=B^\prime\times_BX$ and an induced map
$\phi:|(X^\prime)^.|\rightarrow |X^.|$.
If $\cH\in \cT( |X^.|)\subset \cT(B)$, then we define
$f^*\cH\in \cT(B^\prime)$ as $\phi^*\cH\in \cT(|(X^\prime)^.|)$.
The pull-back of morphisms is defined similarly. In this  way we
obtain a functor $\cT:(orbispaces,representable \:\:maps)\rightarrow groupoids$.

\subsubsection{}

We consider a two-isomorphism
$f\stackrel{\Phi}{\Rightarrow} f^\prime$ between representable maps $f,f^\prime:B^\prime\rightarrow
B$ of orbispaces. If $X\rightarrow B$ is an atlas, and $Y,Y^\prime\rightarrow
B^\prime$ are the atlases obtained by pull-back via $f,f^\prime$, then
$\Phi$ induces a map $\Phi:Y\rightarrow Y^\prime$ which we consider as
a refinement. Note that
 $\phi^\prime\circ |\Phi^.|=\phi:|Y^.|\rightarrow  |X^.|$.
For $\cH\in \cT(|X^.|)\subset \cT(B)$ we define
$\Phi_.(\cH):\phi^*(\cH)\rightarrow |\Phi^.|^*\circ
(\phi^\prime)^*(\cH)$ to be the associated canonical isomorphism,
interpreted as
an isomorphisms $f^*\cH\rightarrow (f^\prime)^*\cH$.

\subsubsection{}

Now we extend the $K$-theory functor (or any other twisted cohomology theory)
to orbispaces.
Let $\cH\in \cT(|X^.|)$ be a twist of $B$ in the sense above.

\begin{ddd}
We define $$K_{Borel}(B,\cH):=K(|X^.|,\cH)\ .$$
\end{ddd}

Let  $f:B^\prime\rightarrow B$ be a map of orbispaces. We use the
notation of \ref{funcrf}. 
\begin{ddd}
We define
$f^*:K_{Borel}(B,\cH)\rightarrow K_{Borel}(B^\prime,f^*\cH)$ to be the map 
$|\phi^.|^*:K(|X^.|,\cH)\rightarrow K(|(X^\prime)^.|,\phi^*\cH)$.
\end{ddd}

Let $\Phi:\cH\rightarrow \cH^\prime$ be a morphism of twists given by
$\phi:u^*\cH\rightarrow (u^\prime)^*\cH^\prime$, where we use the
notation of \ref{tws}.
\begin{ddd}
We define
$\Phi^*: K_{Borel}(B,\cH^\prime)\rightarrow K_{Borel}(B,\cH)$ to be the composition
$$\Phi^*:=(u^*)^{-1}\circ \phi^* \circ (u^\prime)^*\ .$$
\end{ddd} 
Here we us the fact that the refinement map
$u:|Y^.|\rightarrow |X^.|$ is a
homotopy equivalence (see Proposition \ref{heq}), and therefore that $u^*$ is invertible.
We also see that $\Phi^*$ is an isomorphism.

It is straight forward to check that this bi-functor has the required
properties of a twisted cohomology defined on orbispaces as explained
in \ref{rrg1}.

\subsubsection{}

\begin{propo}\label{rt51}
The twisted Borel $K$-theory $K_{Borel}(\dots,\dots)$ is $T$-admissible.
\end{propo}
\proof
We consider the orbispace chart $X:=*\rightarrow [*/\Gamma]$. Then the
corresponding classifying space $|X^.|$ is a countable $CW$-complex of
the homotopy type $B\Gamma$. 
The $T$-duality transformation in $K_{Borel}$ for pairs over
$[*/\Gamma]$ translates to 
the $T$-duality transformation for pairs over $|X^.|$.

In \kpwcite{bunkeschick045} we have shown that the $T$-admissibility of $K$-theory implies
that the $T$-duality transformation is an isomorphism for pairs over bases spaces which are equivalent to finite $CW$-complexes. In fact, this result can be extended to countable complexes as follows.
Let $$W_0\subset W_1\subset \dots \subset W_i\subset\dots \subset W$$ be a filtration of a countable $CW$-complex $W$ by finite sub-complexes.
Let $(E_i,h_i)$, $i=0,1$, be $T$-dual pairs over $W$ and consider twists
$\cH_i\in \cT(E_i)$ such that $[\cH_i]=h_i$.
Let $$T_0:K^{*}(E_0,\cH_0)\rightarrow K^{*-1}(E_1,\cH_1)$$ be the associated $T$-duality transformation. We claim that $T_0$ is an isomorphism of groups.

Let $(E_i(k),h(k))$ be the pairs over $W_k$ obtained by restriction.
We have exact sequences
$$0\rightarrow \lim^1_{\stackrel{k}{\leftarrow}} K^{*-1}(E_i(k),\cH_i(k))\rightarrow
K(E_i,\cH_i)\rightarrow\lim_{\stackrel{k}{\leftarrow}} K^{*}(E_i(k),\cH_i(k))\rightarrow 0$$ for $i=0,1$.
The $T$-duality transformation $T_0$ is compatible with restriction and therefore induces a map of sequences 
$(K^{*}(E_0(k),\cH_0(k)))_{k\ge 0}\stackrel{(T_0(k))_{k\ge 0}}{\rightarrow} (K^{*-1}(E_1(k),\cH_1(k)))_{k\ge 0}
\ .$
Since the complexes $W_k$ are finite, this map is an isomorphism.
We thus obtain a map of short exact sequences
 $$ \begin{array}{cccccc}
0&\rightarrow &\lim^1_{\stackrel{k}{\leftarrow}} K^{*-1}(E_0(k),\cH_0(k))&\longrightarrow&
K(E_0,\cH_0)\\&&(T_0(k))_{k\ge
  0}\downarrow&&T_0\downarrow\\
0&\rightarrow &\lim^1_{\stackrel{k}{\leftarrow}} K^{*-2}(E_1(k),\cH_1(k))&\longrightarrow&
K(E_1,\cH_1)\\[0.7cm]
&&\longrightarrow&\lim_{\stackrel{k}{\leftarrow}}
K^{*}(E_0(k),\cH_0(k))&\rightarrow& 0\\&
&&(T_0(k))_{k\ge 0}\downarrow&&\\
 &&\longrightarrow&\lim_{\stackrel{k}{\leftarrow}}K^{*-1}(E_1(k),\cH_1(k))&\rightarrow& 0\\\end{array}\ .$$
By the five lemma we see that $T_0$ is an isomorphism.
This proves the claim.

We can now apply the claim in order to show that $K_{Borel}$ is
$T$-admissible since the $CW$-complexes $|X^.|$ obtained from
$*\rightarrow  [*/\Gamma]$ for finite groups $\Gamma$ are countable.
\hB

\section{Groupoids and classifying spaces}\label{devo}

\subsection{Continuous cohomology}

\subsubsection{}
We consider a topological groupoid $\cG:\cG^1\Rightarrow \cG^0$ and a topological abelian
group $A$. Then we define a cochain complex of abelian groups
$$\dots \rightarrow C_{cont}^p(\cG,A)\stackrel{\delta}{\rightarrow}
C_{cont}^{p+1}(\cG,A)\rightarrow\dots\ , $$
where 
$$C^0(\cG,A)=C(\cG^0,A)\ ,\quad C_{cont}^p(\cG,A):=C(\underbrace{\cG^1\times
_{\cG^0}\dots \times_{\cG^0}\cG^1}_{p},A)$$ and
\begin{eqnarray*} \lefteqn{(\delta
a)(\gamma_1,\dots,\gamma_{p+1}):=a(\gamma_2,\dots,\gamma_{p+1})}&&\\&&+\sum_{i=1}^p(-1)^i
a(\gamma_1,\dots,\gamma_i\circ\gamma_{i+1},\dots,\gamma_{p+1})+(-1)^{p+1}
a(\gamma_1,\dots,\gamma_p)\ .\end{eqnarray*}
\begin{ddd}
The continuous cohomology $H_{cont}(\cG,A)$ of $\cG$ with values in $A$ is
the cohomology of the complex $(C_{cont}^*(\cG,A),\delta)$.
\end{ddd}
This definition is an immediate extension of the definition of the
continuous cohomology of a topological group.


\subsubsection{}
We now assume that $\cG$ is very proper and {\'e}tale,  and that $A$ admits the structure
of a topological $\R$-vector space.  This allows to multiply $A$-valued continuous functions on a space with $\R$-valued continuous functions.
The following Lemma (see also \cite[Proposition 1]{math.DG/0008064}) generalizes the well-known
fact that the higher cohomology of a finite group with coefficients in a
$\Q$-vector space is trivial. 
\begin{lem}\label{cov}
We have
$H^p(\cG,A)=0$ for $p\ge 1$.
\end{lem}
\proof 
Let $a\in C_{cont}^{p+1}(\cG,A)$ be a cocycle.
We define the continuous cochain $b\in C_{cont}^{p}(\cG,A)$ by
$$b(\gamma_1,\dots,\gamma_p):=(-1)^{p+1}\int_{\cG^{s(\gamma_p)}}
\chi(s(\gamma)) a(\gamma_1,\dots,\gamma_p,\gamma) d\gamma\ ,$$
where $d\gamma$ is the counting measure on
$\cG^{s(\gamma_p)}$, and $\chi$ is the cut-off function given by the very proper condition \ref{verrahgxdhjqw}.
Then by a straight forward computation we have $\delta b=a$.
\hB

\subsection{The Borel construction and $U(1)$-bundles}

\subsubsection{}

We consider a $U(1)$-bundle $E\rightarrow B$ over an orbispace $B$. We choose an orbispace  atlas
$X\rightarrow B$ and get an induced orbispace atlas $Y:=X\times_BE\rightarrow
E$ of $E$. Then we have the groupoids $\cG:X\times_BX\Rightarrow X$ and
$\cE:Y\times_EY\Rightarrow Y$ together with a homomorphism
$\cE\rightarrow \cG$. The latter can be considered as a $U(1)$-bundle
over $\cG$.

It gives rise to a simplicial $U(1)$-bundle
$Y^.\rightarrow X^.$ (using the notation \ref{not}), and thus to an
ordinary $U(1)$-bundle $|Y^.|\rightarrow |X^.|$.

This construction extends in an obvious manner  to a functor
$A_X$ from the category of $U(1)$-bundles over $B$ to $U(1)$-bundles
over $|X^.|$. The morphisms in these categories  here are homotopy classes
of bundle isomorphisms.
The main step in the proof of \ref{main122} is the
following proposition.
\begin{propo}\label{ttz}
$A_X$ is an equivalence of categories.
\end{propo} 
The remainder of the present subsection is devoted to the proof.
It consists of three steps. In the first step we show that $A_X$ is
surjective on the level of sets of isomorphisms classes. Then we show
that it is full. In the last step we show that it is faithful.

\subsubsection{}

We have an equivalence of stacks $B\cong [\cG^0/\cG]$.
Moreover the category of $U(1)$-bundles over $B$ is equivalent to the
category of $U(1)$-bundles over $\cG$.
In fact, given a $U(1)$-bundle $E\rightarrow B$ in stacks we obtain by the construction above
a $U(1)$-bundle $\cE\rightarrow \cG$ in a functorial manner.
In the other direction we funtorially associate to a $U(1)$-bundle
$\cE\rightarrow \cG$ of groupoids a $U(1)$-bundle
$[\cE^0/\cE]\rightarrow [\cG^0/\cG]$ of stacks.

A $U(1)$-bundle $\cE\rightarrow \cG$ in groupoids can equivalently be considered as a $\cG$-equivariant $U(1)$-bundle, i.e. a $U(1)$-bundle
$\cE^0\rightarrow \cG^0$ together with
an action $\cG^1\times_{\cG^0}\cE^0\rightarrow \cE^0$.
Below we will freely switch between these two points of view.

\subsubsection{}
  
If $\cG$ is a topological groupoid then we let $B(\cG)$ denote the
associated simplicial space, and we let $|B(\cG)|$ denote its geometric realization.

In order to prove Proposition \ref{ttz} it suffices to show
that the functor
which associates $|B(\cE)|\rightarrow |B(\cG)|$ to $\cE\rightarrow \cG$
is an equivalence of categories. We will denote it by $A$.

We first show that $A$ induces a surjection on the level of  sets of
isomorphisms classes of objects.

\subsubsection{}

For the following discussion we employ the smooth bundle
$U\rightarrow P \C^\infty$ as a model for the universal
$U(1)$-principal bundle.
To be precise we consider this bundle in the category of $ind$-manifolds such that
$U:=\lim_{\stackrel{n}{\rightarrow}}S^{2n+1}$
and
$P \C^\infty:=\lim_{\stackrel{n}{\rightarrow}}P\C ^n$, and the connecting maps
are in both cases induced by the canonical embeddings $\C^n\rightarrow \C^{n+1}$.

We choose a connection on this $U(1)$ bundle which induces a parallel
transport and a curvature two-form $\omega\in \Omega^2(P \C^\infty)$.
In detail this amounts to choose a compatible family of connections on
the bundles $S^{2n+1}\rightarrow P \C^n$ (e.g. the one induced by the
round metric on the spheres), and the curvature form is interpreted as
a compatible family of two-forms on the family of complex projective spaces, i.e. 
$\omega\in \lim_{\stackrel{n}{\leftarrow}} \Omega^2(P \C^n)$.

A map 
$c:|B(\cG)|\rightarrow   P \C^\infty$
determines a $U(1)$-bundle $c^*U\rightarrow |B(\cG)|$.
Homotopic maps give isomorphic $U(1)$-bundles.  
We want to show that the isomorphism class of $c^*U\rightarrow
|B(\cG)|$ is  in the image of $A$. Let $\bc$ denote the homotopy class
of $c$.

\subsubsection{}

 For all $n\ge 0$
we have a natural map
$$i_n:\Delta^n\times \underbrace{\cG^1\times_{\cG^0}
  \dots\times_{\cG^0} \cG^1}_{n}\rightarrow |B(\cG)|\ .$$
If $(\gamma_1,\dots,\gamma_n)\in \underbrace{\cG^1\times_{\cG^0}
  \dots\times_{\cG^0}\cG^1}_{n}$, then we let
\begin{eqnarray*}\lefteqn{i_n(\gamma_1,\dots,\gamma_n):\Delta^n\cong
\Delta^n\times\{(\gamma_1,\dots,\gamma_n)\}}&&\\&\subset& \Delta^n\times \underbrace{\cG^1\times_{\cG^0}
  \dots\times_{\cG^0}\cG^1}_{n}\stackrel{i_n}{\rightarrow}
|B(\cG)|\ .\end{eqnarray*}

\subsubsection{}
\newcommand{\cir}{{int}}
 We plan to use the parallel transport along one-simplices. 
Furthermore we want to apply
Stokes theorem to the curvature form on three-simplices. Therefore we need a 
representative of $\bc$ which is smooth in the interior
of each simplex. Let $\Delta_\cir^n\subset \Delta^n$ denote the interior of the
standard simplex.
\begin{lem}\label{pars}
The class $\bc$ has a representative 
 $c$ such that for all $n\ge 1$
the composition $c\circ i_n$ induces a continuous map 
$$\underbrace{\cG^1\times_{\cG^0}
  \dots\times_{\cG^0}\cG^1}_{n}\rightarrow C^\infty(\Delta^n_\cir,P
\C^\infty)\ .$$
\end{lem}
\proof
For all $n\ge 1$ we set up one of the usual procedures to smooth out maps 
$\Delta^n\rightarrow P\C^\infty$ in the interior $\Delta^n_\cir\subset
\Delta^n$  without changing the restriction to the boundary.
In this way we obtain  a family of continuous maps 
$C(\Delta^n,P\C^\infty)\rightarrow C^\infty(\Delta^n_\cir,P\C^\infty)\cap C(\Delta^n,P\C^\infty)$. We apply these procedures to the maps
$i_n(\gamma_1,\dots,\gamma_n)$ for all $(\gamma_1,\dots,\gamma_n)\in \underbrace{\cG^1\times_{\cG^0}
  \dots\times_{\cG^0}\cG^1}_{n}$, increasing $n$ from $1$ to $\infty$ inductively. The resulting maps assemble to a representative of $\bc$ with the required properties.
\hB

\subsubsection{}

We define a $U(1)$-bundle $E\rightarrow \cG^0$ by the iterated pull-back
$$\begin{array}{ccccc}
E&\rightarrow&c^*U&\rightarrow &U\\
\downarrow&&\downarrow&&\downarrow\\
\cG^0&\subset&|B(\cG)|&\stackrel{c}{\rightarrow}&P
\C^\infty
\end{array}\ .$$ 
The idea is to define an action of $\cG$ on $E$ so that if we
apply $A$ to the
resulting bundle $\cE\rightarrow \cG$ we get back the isomorphism class of $c^*U\rightarrow |B(\cG)|$.

\subsubsection{}

For $\gamma\in \cG^1$ we have a path
$c\circ i_1(\gamma):\Delta^1\rightarrow P \C^\infty$
from $c(s(\gamma))$ to $c(r(\gamma))$.
We let $\phi(\gamma):E_{s(\gamma)}\rightarrow E_{r(\gamma)}$ denote
the isomorphism such that
$$\begin{array}{ccc}
E_{s(\gamma)}&\stackrel{\phi(\gamma)}{\rightarrow}&E_{r(\gamma)}\\
\|&&\|\\
U_{c(s(\gamma))}&\rightarrow &U_{c(r(\gamma))}\end{array}
\ ,$$
where the lower horizontal arrow is the parallel transport along the
path.

The maps $\phi(\gamma)$, $\gamma\in \cG^1$, combine to a map
$\phi:\cG^1\times _{\cG^0}E\rightarrow E$. This is not yet an action.
In the following we modify this map to make it associative.
In fact, the non-associativity will be measured by a continuous groupoid
cocycle $a$ with coefficients in $U(1)$, and the crucial fact will be that it represents the trivial
cohomology class.

\subsubsection{}
Consider a pair $(\gamma_1,\gamma_2)\in \cG^1\times_{\cG^0}\cG^1$.
We define $$a(\gamma_1,\gamma_2):=\phi(\gamma_1\circ \gamma_2)^{-1}\circ
\phi(\gamma_1)\circ \phi(\gamma_2)\in \Aut(E_{s(\gamma_2)})\cong 
U(1)\ .$$
Note that $a\in C^2_{cont}(\cG,U(1))$ is a cocyle which represents a class
$[a]\in H^2_{cont}(\cG,U(1))$.

\begin{lem}\label{anul}
We have $[a]=0$.
\end{lem}
\proof
We consider the continuous homomorphism $e:\R\rightarrow U(1)$ given by $t\mapsto
\exp(2\pi i t)$. In induces a map of complexes
$e_*:C^2_{cont}(\cG,\R)\rightarrow C^2_{cont}(\cG,U(1))$.
The key to the proof is the observation that the cocycle
$a$ can be lifted to a cocycle $\tilde a\in C_{cont}^2(\cG,\R)$ such
that $e_*\tilde a=a$. By Lemma \ref{cov} we have $[\tilde a]=0$ so
that
$[a]=e_*[\tilde a]=0$, too.

Note that $(\gamma_1,\gamma_2)$ determines a smooth map
$c\circ i_2(\gamma_1,\gamma_2):\Delta^2 \rightarrow  P \C^\infty$.
The restriction of this map  to the boundary of the
simplex determines a piecewise differentiable loop 
in $P\C^\infty$, and $a(\gamma_1,\gamma_2)$ is exactly the holonomy
of the parallel transport along this loop. 
We thus get
$$a(\gamma_1,\gamma_2)=e\left( \int_{\Delta^2}
(c\circ i_2(\gamma_1,\gamma_2))^*\omega\right)\ .$$
We now define the continuous $\R$-valued groupoid-cochain
\begin{equation}\label{adef}\tilde a(\gamma_1,\gamma_2):=\int_{\Delta^2}
(c\circ i_2(\gamma_1,\gamma_2))^*\omega\ .\end{equation}
We claim that $\tilde a$ is a cocycle.
In fact, for $(\gamma_1,\gamma_2,\gamma_3)\in
\cG^1\times_{\cG^0}\cG^1\times_{\cG^0}\cG^1$ the number
$$(\delta \tilde a )(\gamma_1,\gamma_2,\gamma_3)=\tilde a(\gamma_2,\gamma_3) -
\tilde a(\gamma_1\circ \gamma_2, \gamma_3)+
\tilde a(\gamma_1,\gamma_2\circ \gamma_3)-\tilde a(\gamma_1,\gamma_2)$$
is the integral over the boundary of $\Delta^3$ of 
$i_3(\gamma_1,\gamma_2,\gamma_3)^*\omega$. 
 Since $\omega$ is closed, this integral vanishes by Stokes theorem.
\hB 

\subsubsection{}\label{gdr}
By Lemma \ref{anul} we can choose $b\in C_{cont}^1(\cG,U(1))$ such
that \begin{equation}\label{bdef}\delta b=a\ .\end{equation} We now define
$$m(\gamma):=\phi(\gamma) b(\gamma)^{-1}$$
Then it is easy to check that
$m:\cG^1\times_{\cG^0}E \rightarrow E$ is an action.
Let $\cE\rightarrow \cG$ denote the corresponding equivariant $U(1)$-bundle.
\subsubsection{}

Let $F:=|B(\cE)|\rightarrow |B(\cG)|$. 
\begin{lem}\label{tzu}
We have an isomorphism of $U(1)$-bundles $F\cong c^*U$. 
\end{lem}
\proof
 We will prove the assertion by 
explicitly defining an
isomorphism $\psi:F\rightarrow c^*(U)$.

If $(a_0,\dots,a_n)$ are the labels of the vertices of $\Delta^n$, then
let $t_{a_i}$ denote the linear coordinate on $\Delta^n$ which
vanishes at the vertex labeled by $a_i$, and which is equal to $1$ on the
opposite face.

First note that we can find a cochain $\tilde b\in C^1_{Cont}(\cG,\R)$
such that $\delta\tilde b=\tilde a$ and $e(\tilde b)=b$ (using the notation of \ref{gdr}).
Let $\Delta^n$ denote the copy of the standard simplex in $|B(\cG)|$ corresponding to
$$(\gamma_1,\dots,\gamma_n)\in\underbrace{\cG^1\times_{\cG^0}\dots\times_{\cG^0}\cG^1}_{n}\
.$$
The vertices of $\Delta^n$  are naturally labeled by the ordered set
$\{r(\gamma_1),\dots,r(\gamma_n),s(\gamma_n)\}$. 
Let $\Delta_{\circ}^n:=\Delta^n \setminus
\partial_{s(\gamma_n)}\Delta^n$,
where  $\partial_{s(\gamma_n)}\Delta^n$ is the unique face not
containing the vertex labeled by $s(\gamma_n)$.  
We define $\psi$ over the subset $\Delta_\circ^n\times
(\gamma_1,\dots,\gamma_n)\subset |B(\cG)|$ as follows.
By construction the  fiber of $F_{|\Delta_\circ^n\times
(\gamma_1,\dots,\gamma_n)}$ is canonically isomorphic to $E_{s(\gamma_n)}=U_{c(s(\gamma_n))}$.
Each point $s\in \Delta_\circ^n$ can be joined by a linear path
with the vertex with label $s(\gamma_n)$. Let
$\psi(s,(\gamma_1,\dots,\gamma_n)):F_{(s,(\gamma_1,\dots,\gamma_n))}\cong
U_{c(s(\gamma_n))}\rightarrow
U_{c(s,(\gamma_1,\dots,\gamma_n))}$ be given by the parallel
transport along this path multiplied by 
$$e(-t_{s(\gamma_n)}\tilde b(\gamma_n))e(-t_{s(\gamma_n)}t_{s(\gamma_{n-1})}
\tilde b(\gamma_{n-1}))\dots e(-t_{s(\gamma_n)}\dots t_{s(\gamma_1)}
\tilde b(\gamma_1))\ .$$  
We use the construction for all $n\ge 1$ and points $(\gamma_1,\dots,\gamma_n)\in\underbrace{\cG^1\times_{\cG^0}\dots\times_{\cG^0}\cG^1}_{n}$. 
It is now easy to check that $\psi$ is an everywhere  defined
continuous bundle isomorphism. \hB

 This finishes the proof of the fact that $A$ is surjective on the level of sets of isomorphism classes of objects.

\subsubsection{}\label{rtd1}

Our next task is to show that $A$ is full.
We consider the following intermediate construction.
Let $\cE\rightarrow \cG$ be a $U(1)$-bundle.
Then we have a cartesian diagram
\begin{equation}\label{xdx}\begin{array}{cccccc}
|B(\cE)|&\stackrel{\cong}{\rightarrow}&c^*U&\rightarrow &U\\
\downarrow&&\downarrow&&\downarrow\\
|B(\cG)|&\stackrel{=}{\rightarrow} & |B(\cG)|&\stackrel{c}{\rightarrow}&P
\C^\infty
\end{array}\ ,\end{equation}
where $c$ is uniquely determined up to homotopy. 
After a further homotopy we can assume
that $c$ satisfies the condition of Lemma \ref{pars}.
We apply to this map $c$ the construction of the first part of the
proof and obtain a 
$U(1)$-bundle $\tilde \cE\rightarrow \cG$.

\subsubsection{}
\begin{lem}\label{zut}
We have $\tilde \cE\cong \cE$ as $U(1)$-bundles over $\cG$.
\end{lem}
\proof
Let $E,\tilde E\rightarrow \cG^0$ be the underlying $U(1)$-bundles.
Note that (\ref{xdx}) induces a canonical isomorphism
$\Psi:\tilde E \stackrel{\sim}{\rightarrow} E$ as $U(1)$-principal bundles over
$\cG^0$. We must compare the action $\tilde m$ of $\cG$ on $\tilde E$ with the
original action $m$ on $E$. The difference between these two  actions
is measured by the continuous cocycle $h\in C^1_{cont}(\cG,U(1))$
defined by
$$h(\gamma)=\Psi^{-1}\circ m(\gamma)^{-1}\circ \Psi \circ  \tilde m(\gamma)\in \Aut(
\tilde E_{s(\gamma)})\cong U(1)\ .$$
The cohomology class of this cocycle is the obstruction against making $\Psi$ equivariant
by multiplying it by a $U(1)$-valued function on $\cG^0$.
 
\subsubsection{}
\begin{lem}\label{hnul}
We have $[h]=0$.
\end{lem}
\proof
The key  is again the construction of a lift of $h$ to a cocycle
$\tilde h\in C_{cont}^1(\cG,\R)$ such that $e_*(\tilde h)=h$.
By Lemma \ref{cov} we then have $[h]=e_*([\tilde h])=0$. 

We consider $\gamma\in \cG^1$.
 It induces a smooth path $c\circ i_1(\gamma):\Delta^1 \rightarrow  P\C^\infty$ and
therefore
a parallel transport $\phi(\gamma):U_{c(s(\gamma))}\rightarrow
U_{c(r(\gamma))}$. We have $\tilde m(\gamma)=\phi(\gamma)
b(\gamma)^{-1}$,
where $b$ is as in (\ref{bdef}).
As in the proof of Lemma \ref{tzu}  will again use the cochain $\tilde b\in
C_{cont}(\cG,\R)$ such that $\delta \tilde b=\tilde a$ and
$b=e_*(\tilde b)$.

The identification $|B(\cE)|\cong c^*U$ induces a trivialization
$i_1(\gamma)^*U\cong \Delta^1\times E_{s(\gamma)}$.  
If $\alpha(\gamma)$ denotes the connection-one form in this
trivialization, then we can write
$$\phi(\gamma)=e\left( \int_{\Delta^1} \alpha(\gamma)\right)\ .$$
By construction we have
$h(\gamma)=e\left( \int_{\Delta^1} \alpha(\gamma)\right)b(\gamma)^{-1}$.
We define the cochain $\tilde h\in C_{cont}^1(\cG,\R)$
$$\tilde h(\gamma):=\int_{\Delta^1} \alpha(\gamma)-\tilde b(\gamma)\ .$$
It satisfies $e_*(\tilde h)=h$.
We claim that $\tilde h$ is in fact a cocycle. 
Let $(\gamma_1,\gamma_2)\in \cG^1\times_{\cG^0}\cG^1$. The identification $|B(\cE)|\cong c^*U$ induces a trivialization
$(c\circ i_2(\gamma_1,\gamma_2))^*U\cong \Delta^2\times U_{E_{s(\gamma_2)}}$.
Let $\alpha(\gamma_1,\gamma_2)$ denote the connection one-form in this
trivialization. 
Then we have
$$\delta \tilde h(\gamma_1,\gamma_2)=\int_{\partial
  \Delta^2}\alpha(\gamma_1,\gamma_2) - \delta\tilde
b(\gamma_1,\gamma_2)\ .$$
By Stoke's theorem the first term of the right-hand side
 is
equal to
$$\int_{\Delta^2} d\alpha(\gamma_1,\gamma_2)\ .$$
Now the claim follows in view of
$ d\alpha(\gamma_1,\gamma_2)=(c\circ i_2(\gamma_1,\gamma_2))^*\omega$, $\delta\tilde b=\tilde a$,
and  (\ref{adef}). \hB

\subsubsection{}
By Lemma \ref{hnul} we can choose a cochain
$f\in C^0_{cont}(\cG,U(1))$ such that $\delta f=h$.
If we define the isomorphism $\tilde \Psi:\tilde E\rightarrow E$ by
$\tilde \Psi(x)=\Psi(x) f^{-1}(x)$ then
$\tilde \Psi$ is  $\cG$-equivariant.

\subsubsection{}
We now finish the proof of the fact that $A$ is full. To this end 
we consider $U(1)$-bundles
$\cE,\cE^\prime\rightarrow \cG$ and an isomorphism of $U(1)$-bundles
$\Lambda:|B(\cE^\prime)|\rightarrow |B(\cE)|$ over $|B(\cG)|$.
We must show that $\Lambda$ can be written as
$A(\lambda)$ for some $ \lambda:\cE^\prime\rightarrow \cE$ over $\cG$.
We apply to $\cE$ and $\cE^\prime$ the intermediate construction started in \ref{rtd1}, where we use the same map $c:|B(\cG)|\rightarrow P \C^\infty$ in both cases. We
obtain a chain of isomorphisms
$$\cE\stackrel{\tilde \Psi}{\cong}\tilde \cE=
 \tilde
\cE^\prime\stackrel{\tilde \Psi^\prime}{\cong}\cE^\prime\ .$$
Let $\cE\stackrel{\lambda}{\cong}\cE^\prime$ be the composition. 

In general $A(\lambda)$ is not equal to $\Lambda$ (recall that we consider homotopy classes). But the following result shows that we can find an automorphism $\phi$ of $\cE$ such that
$A(\lambda\circ \phi)=\Lambda$.

\subsubsection{}

Let $\phi:\cG^0\rightarrow U(1)$ be a $\cG^1$-invariant function.
We can interpret $\phi$ as an automorphism of the
$U(1)$-bundle $\cE\rightarrow \cG$.
Applying the classifying space functor we get an automorphism
$|B(\phi)|$ of the $U(1)$-bundle $|B(\cE)|\rightarrow |B(\cG)|$,
i.e. a function $|B(\phi)|:|B(\cG)|\rightarrow U(1)$.

\begin{lem}\label{suro}
Every homotopy class of maps $[|B(\cG)|,U(1)]$ has a
representative of the form $|B(\phi)|$ for a $\cG^1$-invariant function
$\phi:\cG^0\rightarrow U(1)$.
\end{lem}
\proof
We consider a homotopy class of maps $|B(\cG)|\rightarrow U(1)$ and choose a representative $\tilde f$.
The restriction of $\tilde f:|B(\cG)|\rightarrow U(1)$ to $\cG^0\subset |B(\cG)|$ gives a function
$\tilde \phi:\cG^0\rightarrow U(1)$. In general it is not
$\cG^1$-invariant.

We consider $\tilde \phi\in C_{cont}^0(\cG,U(1))$. Then the non-invariance
is measured by $h:=\delta \tilde \phi\in C_{cont}^1(\cG,U(1))$.

We have $h(\gamma)=\phi(r(\gamma))\phi(s(\gamma))^{-1}$.
We now construct a lift $\tilde h\in C_{cont}^1(\cG,\R)$
as follows. Let $\gamma\in \cG^1$. It gives rise to a path
$i_1(\gamma):\Delta^1\rightarrow |B(\cG)|$.
The restriction $i_1(\gamma)^*\tilde f$ has a lift to an $\R$-valued
function $\kappa(\gamma):\Delta^1\rightarrow \R$. The difference  $\tilde h(\gamma):=\kappa(\gamma)(1)-\kappa(\gamma)(0)$ is independent of the choice of the lift.

We claim that $\delta\tilde h=0$. This follows from the fact that
$\tilde f$ is defined on the image of
$i_2(\gamma_1,\gamma_2):\Delta^2\rightarrow |B(\cG)|$ for all
composeable $\gamma_1,\gamma_2\in \cG^1$.
By Lemma \ref{cov} we can find a function $a\in C_{cont}^0(\cG,\R)$
such that $\delta a=\tilde h$.
We now define the $\cG^1$-invariant $U(1)$-valued function
$$\phi=\tilde \phi \exp(-2\pi i a)\ .$$
We can consider $a$ as an $\R$-valued continuous function defined on
the closed subset $\cG^0\subset  |B(\cG)|$. Let $\tilde
a:|B(\cG)|\rightarrow \R$ be any continuous extension, and set
$f:=\tilde f\exp(-2\pi i \tilde a)$. Then clearly $[f]=[\tilde f]$.
It remains to show that $[f]=[|B(\phi)|]$.

Note that
$i_n(\gamma_1,\dots,\gamma_n)^*B(\phi)=\phi(s(\gamma_n))=\phi(r(\gamma_i))$
for all $i=1,\dots,n$.
We now consider the function $g:|B(\cG)|\rightarrow U(1)$ defined by
$g=f B(\phi)^{-1}$. It has the property that $g_{|\cG^0}=1$.
We must show that $g$ is homotopic to the constant function,
or equivalently, that it admits a lift to an $\R$-valued function.
In fact, in this case $[f]=[|B(\phi)|]$.

We have a natural map $p:|B(\cG)|\rightarrow \cG^0/\cG^1$
(the target is the quotient space of $\cG^0$ with respect to the equivalence relation generated by $\cG^1$)
given by $p(\sigma,(\gamma_1,\dots,\gamma_n)):=s(\gamma_n)$,
where $\sigma\in \Delta^n$. The fibre of $p$ over the class
$[x]\in \cG^0/\cG^1$ is homotopy equivalent to the classifying
space $|B(\cG^x_x)|$.
Since $\cG^x_x$ is a finite group we have
$H^1(|B(\cG^x_x)|,\Z)=0$.
This shows that the restriction of the $U(1)$-valued function $g$ to $p^{-1}([x])$ admits a lift
to an $\R$-valued function which is unique up to an additive integer.

Let $[x] \in \cG^0/\cG^1$ and $\gamma\in \cG^1$ such that $s(\gamma)\in
[x]$. Let
$\tilde g_{[x]}$ be a lift of
$g_{|p^{-1}([x])}$. Then we have $\tilde g_{[x]}(r(\gamma))-\tilde g_{[x]}(s(\gamma))=
\kappa(\gamma)(1) -\kappa(\gamma)(0)-\tilde a(r(\gamma))+\tilde
a(s(\gamma))=
\tilde h(\gamma)-\tilde a(r(\gamma))+\tilde
a(s(\gamma))=0$.
This allows us to normalize the lift $\tilde g_{[x]}$ such that
$(\tilde g_{[x]})_{|[x]}=0$. These normalized lifts fit together
to a lift $\tilde g:|B(\cG)|\rightarrow \R$ of $g$. \hB

This finishes the proof of the fact that $A$ is full.
Note that this implies that $A$ is injective on the level of sets of isomorphism classes of objects.
\subsubsection{}

In the final step of the proof of Proposition \ref{ttz} we show that $A$ is faithful. It suffices to show that $A$ is injective on the group of automorphisms of a $U(1)$-bundle $\cE\rightarrow \cG$. Via  a mapping torus construction we can translate this assertion to the injectivity of $A$ on the set of isomorphism classes of $U(1)$-bundles over $S^1\times \cG$. Therefore faithfulness is implied by the preceeding results. This finishes the proof of Proposition \ref{ttz}. \hB

\subsection{The Borel construction for pairs}

\subsubsection{}

In this subsection we finish the proof of Theorem \ref{main122}.
Let $Y\rightarrow B$ be an atlas of an orbispace $B$.
Recall that
$PA_Y:P(B)\rightarrow P(|Y^.|)$ maps the pair
$(E,h)$ to $(|X^.|,h)$, where $X:=E\times_BY\rightarrow E$ is the induced atlas
of $E$, $|X^.|\rightarrow |Y^.|$ is the induced $U(1)$-principal bundle, and
$h\in H^3(|X^.|,\Z)\cong H^3(E,\Z)$.

We must show that $PA_Y$ induces an isomorphism on the level of isomorphism classes pairs. Since the construction is functorial it is clear that $PA_Y$ descends to isomorphism classes.

We first show that it is surjective. Consider a pair $(F,h)$ over $|Y^.|$.
Then by Proposition \ref{ttz} we find a $U(1)$-bundle $E\rightarrow B$ such that $|X^.|\cong F$ as $U(1)$-bundles over $|Y^.|$. Using this isomorphism we consider $h\in H^3(E,\Z)$. It follows that
$A_Y$ maps $(E,h)$ to $(F,h)$. Hence, $PA_Y$ hits all isomorphism classes.

We now consider two pairs $(E_i,h_i)$, $i=0,1$ over $B$. We assume that they become isomorphic under $PA_Y$, i.e. we have an isomorphism of $U(1)$-bundles
$\phi:|X_0^.|\rightarrow |X_1^.|$ such that $\phi^*h_1=h_0$.
We apply again Proposition \ref{ttz} in order to find
an isomorphism $\Phi:E_0\rightarrow E_1$ such that $PA_Y(\Phi)$ is homotopic to $\phi$. It therefore gives an isomorphism of pairs
$(E_0,h_0)\cong (E_1,h_1)$.
This shows that $PA_Y$ is injective.\hB

\section{Examples}\label{exex}
\subsection{$\Gamma$-Points - cyclic groups}\label{gp1}

\subsubsection
Let $\Gamma$ be a finite group. Let $\Gamma$ act on the one point
space $*$ and consider the orbispace $B:=[*/\Gamma]$. The map
$*\mapsto [*/\Gamma]$ is an atlas. The associated groupoid is $\cG:\Gamma\Rightarrow *$, and $B(\cG)$ is the usual bar construction on $\Gamma$. We have $|B(\cG)|\cong B\Gamma$.

\subsubsection{}

 The group of characters of $\Gamma$ can be identified with the group cohomology  $H^1(\Gamma,U(1))$.
Let $\chi\in H^1(\Gamma,U(1))$ be a character.
It induces an action of $\Gamma$ on $U(1)$. We obtain a
$U(1)$-principal bundle $E:=[U(1)/\Gamma]\rightarrow B$.
In order to extend $E$ to a pair over $B$ we must choose a class
$h\in H^3(E,\Z)$. We use the Gysin sequence in order to get some information about this cohomology group.

\subsubsection{}

The topology of the bundle  $E\rightarrow B$ enters into the Gysin sequence through its first Chern class. In order to describe this class in terms of the character $\chi$ we consider
the boundary operator  of the long exact sequence in group cohomology 
 associated to the sequence of coefficients
$$0\rightarrow \Z\rightarrow \R\rightarrow U(1)\rightarrow 0\ .$$ It
provides an isomorphism
$$\delta:H^1(\Gamma,U(1))\stackrel{\sim}{\rightarrow}
H^2(\Gamma,\Z)\cong H^2(B\Gamma,\Z)\cong H^2(B,\Z)\ .$$
Let $c_1(E)\in H^2(B,\Z)$ denote the first Chern class of $E$.
We then have
$$c_1(E)=\delta(\chi)\ .$$

\subsubsection{}

Since $\Gamma$ is finite we have $H^1(B\Gamma,\Z)=H^1(B,\Z)=0$.
The relevant part of the Gysin sequence has the form
$$0\rightarrow H^3(B,\Z)\stackrel{\pi^*}{\rightarrow} H^3(E,\Z)\stackrel{\pi_!}{\rightarrow} H^2(B,\Z)\stackrel{\dots\cup c_1(E)}{\rightarrow} H^4(B,\Z)\rightarrow \dots\ .$$

\subsubsection{}\label{zwde}
Let us from now on assume that $\Gamma$ is the cyclic group $\Z/n\Z$.
We identify $\hat \Gamma\cong \Z/n\Z$ such that the character
corresponding to $[q]\in \Z/n\Z$ is given by
$$\chi([p])=\exp(\frac{2\pi i p q}{n})\ .$$  

The cohomology of $B\Gamma$ is given by
$$\begin{array}{|c|c|}\hline
i&H^i(B\Gamma,\Z)\\\hline
0&\Z\\\hline
2l-1&0\\\hline
2l&\Z/n\Z\\\hline
\end{array}\ ,$$
where $l\ge 1$.

Under this identification we have
$c_1(E)=[q]$. The Gysin sequence specializes to
$$0\rightarrow H^3(E,\Z)\stackrel{\pi_!}{\rightarrow}
\Z/n\Z\stackrel{[q]}{\rightarrow}
\Z/n\Z\rightarrow \dots$$
so that $$H^3(E,\Z)\cong \{[s]\in \Z/n\Z\:|\: n|sq\}\subset \Z/n\Z \ .$$
We fix a class $h=[s]$ in this group.

\subsubsection{} 

We can now calculate the $T$-dual pair $(\hat E,\hat h)$.
Note that by \kpwcite{bunkeschick045}, Lemma 2.12, we have 
$c_1(\hat E)=-\pi_!(h)$. Therefore, we have
$c_1(\hat E)=[-s]\in \Z/n\Z\cong H^2(B,\Z)$. 
We can determine $\hat h$ by the condition $\hat \pi_!(\hat h)=-c_1(E)$.
The relevant part of the Gysin sequence for $\hat E$ has the form
$$0\rightarrow H^3(\hat E,\Z)\stackrel{\hat \pi_!}{\rightarrow}
\Z/n\Z\stackrel{[-s]}{\rightarrow}
\Z/n\Z\rightarrow \dots \ ,$$ so that
$$H^3(\hat E,\Z)=\{[r]\in \Z/n\Z\:|\: n|sr\}\subset \Z/n\Z\ ,$$
and we have
$\hat h = [-q]$.
 
\subsubsection{}

Note that the stack $E=[U(1)/\Z/n\Z]$ is equivalent to a space which
is homeomorphic to $U(1)$. But the action of $U(1)$ on this space is not free.
Let us assume that $(q,n)=1$. Then we have $H^3(E,\Z)=0$ and thus $h=0$.
The dual bundle is then given by the orbispace $\hat E=[U(1)/\Z/n\Z]$, where
the group $\Z/n\Z$ now acts trivially. This orbispace is not equivalent to a space. We have $H^3(\hat E,\Z)\cong \Z/n\Z$, and $h=[-q]$.
This example shows that in general the $T$-dual of a space with a non-free $U(1)$-action is an orbispace which is not equivalent to a space anymore.

\subsubsection{}

We now calculate the twisted Borel $K$-groups for $E$ and $\hat E$.
As predicted by the general theory they turn out to be isomorphic (up to degree-shift).
We keep the assumption $(n,q)=1$.

Since $h=0$ and $E\cong U(1)$ we have 
$$\begin{array}{|c|c|}\hline
i&K^i_{Borel}(E,\cH)\\\hline
2l-1&\Z\\\hline
2l&\Z\\\hline
\end{array}\ ,$$
where $l\in \Z$ and $\cH$ is a trivializable twist.
\newcommand{\ori}{{\tt or}}

\subsubsection{}

We employ the Mayer-Vietoris sequence in order to calculate
$K_{Borel}^*(\hat E,\hat \cH)$, where $\hat \cH$ is a twist of $\hat E\cong U(1)\times [*/\Z/n\Z]$ classified by $\hat h$. We fix the atlas $*\rightarrow [*/\Z/n\Z]$. Then $X:=U(1)\times *\rightarrow U(1)\times [*/\Z/n\Z]$ is an atlas of $\hat E$. We get $|X^.|\cong U(1)\times B\Z/n\Z$. We have
$\hat h=\ori_{U(1)}\times [-q]$, where
$\ori_{U(1)}\in H^1(U(1),\Z)$ is the positive generator, and
$[-q]\in H^2(B\Z/n\Z,\Z)\cong \Z/n\Z$. We can assume that $\hat \cH$ is a twist on $|X^.|$.
We decompose $U(1)$ into the union of an upper and a lower hemisphere $I^\pm$.
The restriction of $\hat \cH$ to $I^\pm\times B\Z/n\Z$ is trivializable.

\subsubsection{}

We have a ring isomorphism $K(B\Z/n\Z)\cong R(\Z/n\Z)_{(I)}$,
where $I\subset R(\Z/n\Z)$ is the dimension ideal in the representation ring of $\Z/n\Z$, and $(\dots)_{(I)}$
denotes the  $I$-adic completion.
In particular we have $K^1(B\Z/n\Z)\cong \{0\}$.
 We have a natural map
$\Z/n\Z\rightarrow K(B\Z/n\Z)$ which associates to
$[q]$ the class of the line bundle over $B\Z/n\Z$  associated to the character
$[s]\mapsto \exp(2\pi i\frac{sq}{n})$.

\subsubsection{}

We can now write out the 
Mayer-Vietoris sequence in twisted $K$-theory associated to the decomposition
$$|X^.|\cong (I^+\times B\Z/n\Z)\cup (I^-\times B\Z/n\Z)\ .$$
\begin{eqnarray*}\lefteqn{0\rightarrow  K^0_{Borel}(\hat E,\hat \cH)\rightarrow}&&\\&&
K^0(B\Z/n\Z)\oplus K^0(B\Z/n\Z)\stackrel{\left(\begin{array}{cc}1&1\\-[-q]&-1\end{array}\right)}{\rightarrow}K^0(B\Z/n\Z)\oplus K^0(B\Z/n\Z)
\\&&\rightarrow
K^1_{Borel}(\hat E,\hat \cH)\rightarrow 0\ .\end{eqnarray*}
Here, since $I^\pm$ is contractible and the restriction of the twist is trivializable, we  identify $K(I^\pm\times B\Z/n\Z,\hat \cH_{|I^\pm\times B\Z/n\Z})$
with $K(B\Z/n\Z)$. The appearance of $[-q]$ instead of $-1$ in the lower left corner of  the matrix is due to the presence of twists.
We now use the isomorphism
$K(B\Z/n\Z)\cong R(\Z/n\Z)_{(I)}$ and calculate that
$$K^0_{Borel}(\hat E,\hat \cH)\cong \ker(([-q]-1):R(\Z/n\Z)_{(I)}\rightarrow R(\Z/n\Z)_{(I)})\cong \Z$$ and
$$K^1_{Borel}(\hat E,\hat \cH)\cong \coker(([-q]-1):R(\Z/n\Z)_{(I)}\rightarrow R(\Z/n\Z)_{(I)})\cong \Z\ .$$
Therefore we get
 $$\begin{array}{|c|c|}\hline
i&K^i_{Borel}(\hat E,\hat \cH)\\\hline
2l-1&\Z\\\hline
2l&\Z\\\hline
\end{array}$$
as predicted by the $T$-duality isomorphism.

\subsection{Seifert fibrations}

\subsubsection{}

In this subsection we consider $T$-duality of $U(1)$-bundles
over certain two-dimensional orbispaces. In order to describe such an orbispace $B$ we  
 fix numbers $r,g\in \nat_0$, and an element
$(n_1,\dots,n_r)\in (\Z\setminus \{0\})^r$. 
We set $n_0:=1$.
We consider $\Gamma_i:=
\Z/n_i\Z$ as a subgroup of $U(1)$ via
$[q]\mapsto \exp(2\pi i\frac{q}{n_i})$.

Let $\Sigma$ be an oriented surface of genus $g$.
We fix pairwise distinct points $p_0,p_1,\dots,p_r\in
\Sigma$. We further choose orientation preserving identifications
$(\bar U_i,p_i)\cong (D^2,0)$ of suitable pairwise disjoint closed pointed neighborhoods $\bar U_i$ of $p_i$ for all
$i=0,\dots r$. 
The group
 $\Gamma_i$ acts  naturally on the disk $\tilde D\subset \C$.
We consider the associated branched covering $\tilde D\rightarrow D$, $z\mapsto z^{|n_i|}$, 
and let $\tilde{\bar U}_i\rightarrow \bar U_i$ be the branched covering induced via our identification $\bar U_i\cong D$.

\subsubsection{}
This data determines a topological groupoid $\cG$ which represents the orbispace $B:=[\cG^1/\cG^0]$.
 Let $\Sigma^0:=\Sigma\setminus \bigcup_{i=0}^r U_i$, where $U_i\subset \bar U_i$ denotes the interior.
We define
$$\cG^0:=\Sigma^0  \sqcup \bigsqcup_{i=0}^r
\tilde{\bar U}_i\ .$$
The set of morphisms is defined as follows.
First of all the restriction of $\cG$ to
$\Sigma^0$ is the trivial groupoid.
The restriction of $\cG$ to $\tilde{\bar U}_i$ is the action groupoid
of the $\Gamma_i$-action on $\tilde{\bar U}_i$, i.e.
$\Gamma_i\times  \tilde{\bar U}_i\Rightarrow \tilde{\bar U}_i$.
It remains to describe the morphisms over the overlaps.
A point $s^\Sigma\in \partial \Sigma^0$ determines an index $i$ and a point
$\bar s\in \bar U_i$. For any lift $\tilde s\in \tilde{\bar U}_i$ of
$\bar s$ we require that
there is exactly one morphism
$s^\Sigma\rightarrow \tilde s$ in $\cG^1$.
As a topological space $\cG^1$ is fixed by the requirement that
$s:s^{-1}(\partial \Sigma^0)\rightarrow \partial \Sigma^0$ is a connected covering over each connected component of $\partial \Sigma^0$, where $s:\cG^1\rightarrow \cG^0$ is the source map.

In fact, this groupoid describes an orbispace structure
on $\Sigma$ with singular points $p_1,\dots,p_r$ of multiplicity
$n_1,\dots,n_r$. The point $p_0$ will be used later in order to
introduce a non-trivial topology on $U(1)$-bundles over $B$ in the case $r=0$.

\subsubsection{}\label{con}
We now describe $U(1)$-bundles over $B$. To this end we
 choose a number $c\in \Z$ and an element
$(\chi_1,\dots,\chi_r)\in \hat \Gamma_1\times\dots\times \hat \Gamma_r$.
This data together with additional choices (the $\phi_i$ introduced below) determines a $U(1)$-bundle $E\rightarrow B$ as follows.
We will describe it as a quotient $E:=[\cE/\cG^1]$, where $\cE\rightarrow \cG$ is an equivariant $U(1)$-bundle. It is given by a $U(1)$-bundle
$\cE\rightarrow \cG^0$ together with an action $\cG^1\times_{\cG^0}\cE\rightarrow \cE$.
We set $\cE:=U(1)\times \cG^0$. The data fixed  above determines the action of
$\cG$. On $\cE_{|\tilde {\bar U}_i}$ we let $\Gamma_i$ act on the fibre
with character $\chi_i$. 

For all $i=1,\dots,r$ we choose a map
$\phi_i:\partial \tilde{\bar U}_i\rightarrow U(1)$ such that
$\phi_i(\gamma \tilde s)=\chi_i(\gamma) \phi_i(\tilde s)$, $\gamma\in \Gamma_i$. We identify $\hat \Gamma_i\cong \Z/n_i\Z$ such that
$[q]\in \Z/n_i\Z$ corresponds to the character $[p]\mapsto \exp(2\pi i\frac{pq}{n_i})$.
Note that in $\hat \Gamma_i\cong \Z/n_i\Z$ we have $[\deg(\phi_i)]= \chi_i$.
Here in order to define the degree $\deg(\phi_i)\in \Z$, we choose the orientation
of $\partial \tilde {\bar U}_i$ as the boundary of the oriented disk
$ \tilde {\bar U}_i$.
  Furthermore note that two  choices of $\phi_i$  differ by a function
$\partial \bar U_i\rightarrow U(1)$. Thus we can realize all  elements of the residue class of $\chi$ as $\deg(\phi_i)$ for an appropriate choice of $\phi_i$.

We let the  morphism $s^\Sigma\rightarrow \tilde s$ act as
multiplication by $\phi_i(\tilde s)$, if $s^\Sigma$ is in the $i$th component
of $\partial \Sigma^0$, $i=1,\dots,r$.

Finally, we take a function $u:\partial \bar U_0\rightarrow U(1)$ of degree
$c$ and let the morphism $s^\Sigma\rightarrow s$ act by
multiplication by $u(s)$, if $s^\Sigma$ is in the zero-component of $\partial \Sigma^0$.

\subsubsection{}

If $\chi_i$ are generators of $\hat \Gamma_i$ for all $i=1,\dots r$, then
$E$ is a space. Otherwise  $E$ is an orbispace which is not equivalent
to a space.

\subsubsection{}

We first compute $H^*(B,\Z)$ using a Mayer-Vietoris sequence.
We obtain
\begin{eqnarray*}\dots\rightarrow \bigoplus_{i=0}^r H^{*-1}(\partial \bar U_i,\Z)\rightarrow
H^*(B,\Z)\rightarrow H^*(\Sigma^0,\Z)\oplus\bigoplus_{i=0}^r
H^*(B\Gamma_i,\Z)&&\\\rightarrow  \bigoplus_{i=0}^r H^{*}(\partial \bar U_i,\Z)\rightarrow\dots&&\
.\end{eqnarray*}
We have a canonical identification $H^2(B\Gamma_i,\Z)\cong \hat \Gamma_i$.
The fixed embedding $\Gamma_i\hookrightarrow U(1)$ induces a map
$$B\Gamma_i\rightarrow BU(1)\cong K(\Z,2)$$ and therefore
a generator $c_i\in H^2(B\Gamma_i,\Z)$.
The multiplication with the powers of $c_i$  provides the isomorphisms $\hat \Gamma_i\cong H^{2l}(B\Gamma_i,\Z)$. Furthermore, $H^{2l-1}(B\Gamma_i,\Z)\cong\{0\}$.

\subsubsection{}\label{h2b}
The Mayer-Vietoris sequence now gives the following information.
$$\begin{array}{|c|c|}\hline
l&H^l(B,\Z)\\\hline
0&\Z\\\hline
1&\Z^{2g}\\\hline
2&0\rightarrow
\Z \stackrel{\delta}{\rightarrow}
H^2(B,\Z)\rightarrow \bigoplus_{i=1}^r \hat
\Gamma_i\rightarrow 0\\\hline
2l+1,l\ge 1&0\\\hline
2l,l\ge 2&\bigoplus_{i=1}^r \hat
\Gamma_i\\
\hline
\end{array}\ .$$

The data chosen in the construction \ref{con} provides a split $s$ of the exact
sequence for $H^2(B,\Z)$. In fact, given
$(\chi_1,\dots,\chi_r)\in \bigoplus_{i=1}^r \hat
\Gamma_i$ we construct the line bundle $E\rightarrow B$ associated to this tuple
and $c=0$. Then we set $s(\chi_1,\dots,\chi_r):=c_1(E)$.
It will follow from the calculations in \ref{c1232} that this gives a split.
Since there is no non-trivial homomorphism from
$ \bigoplus_{i=1}^r \hat
\Gamma_i$ to $\Z$ the split $s$ is independent of the choices.
Therefore we can unambiguously write 
$$H^2(B,\Z)\cong \Z\oplus \bigoplus_{i=1}^r \hat
\Gamma_i\ .$$
We will write elements in the form
$(e,(\kappa_1,\dots,\kappa_r))$.

\subsubsection{}\label{c1232}

By Proposition \ref{ttz} the topological type of the $U(1)$-bundle $E\rightarrow B$ is classified by its first Chern class $c_1(E)$. In the following paragraph we calculate this invariant. To this end we consider the following part of the Gysin sequence of $\pi:E\rightarrow B$:
$$\Z\cong H^0(B,\Z)\stackrel{c_1(E)}{\rightarrow} H^2(B,\Z)\stackrel{\pi^*}{\rightarrow} H^2(E,\Z)\ .$$ 
We see that we can calculate $c_1(E)$  by determining the corresponding generator of the kernel of $\pi^*:H^2(B,\Z)\rightarrow H^2(E,\Z)$.

We obtain information on  $H^2(E,\Z)$ using the Mayer-Vietoris sequence. The relevant part has the form
\begin{eqnarray*}
H^1(U(1)\times \Sigma^0,\Z)\oplus \bigoplus_{i=0}^rH^1([U(1)/_{\chi_i}\Gamma_i],\Z)\stackrel{\alpha}{\rightarrow} \bigoplus_{i=0}^r  H^1(U(1)\times \partial \bar U_i,\Z)&&\\\rightarrow H^2(E,\Z)\rightarrow&&\\
H^2(U(1)\times \Sigma^0,\Z)\oplus \bigoplus_{i=0}^rH^2([U(1)/_{\chi_i}\Gamma_i],\Z)\stackrel{\beta}{\rightarrow} \bigoplus_{i=0}^r  H^2(U(1)\times \partial \bar U_i,\Z)\ .
\end{eqnarray*}
The known cohomology groups are
\begin{eqnarray*}
 H^1(U(1)\times \Sigma^0,\Z)&\cong&1_{U(1)}\times H^1(\Sigma^0,\Z)\oplus \ori_{U(1)}\times (1_{\Sigma^0})\Z\\
H^1(\Sigma_0,\Z)&\cong&\Z^{2g+r}\\
 H^1(U(1)\times \partial \bar U_i,\Z)&\cong&(1_{U(1)}\times\ori_{\partial \bar U_i})\Z\oplus (\ori_{U(1)}\times 1_{\partial \bar U_i})\Z\\
H^1([U(1)/_{\chi_i}\Gamma_i],\Z)&\cong& \Z\\
H^2(U(1)\times \Sigma^0,\Z)&\cong&\ori_{U(1)}\times H^1(\Sigma^0,\Z)\\
H^2(U(1)\times \partial \bar U_i,\Z)&\cong&(\ori_{U(1)}\times \ori_{\partial \bar U_i})\Z\\
H^2([U(1)/_{\chi_i}\Gamma_i],\Z)&\cong&\hat \Gamma_i/\chi_i \hat \Gamma_i\ ,
\end{eqnarray*}
where the definition of $\hat \Gamma_i/\chi \hat \Gamma_i$ uses the ring structure on $\hat \Gamma_i\cong \Z/n\Z$.

The map $\beta$ vanishes on the torsion subgroups $H^2( [U(1)/_{\chi_i}\Gamma_i],\Z) $. The range of the restriction of $\beta$ to $H^2(U(1)\times \Sigma^0,\Z)$ has rank $r$. We see that
$$\ker(\beta)\cong \Z^{2g}\oplus \bigoplus_{r=1}^r \hat \Gamma_i/_{\chi_i} \hat \Gamma_i\ .$$

We now determine the cokernel of $\alpha$.
We proceed in stages. We first determine the cokernel of the restriction of $\alpha$ to $1_{U(1)}\times H^1(\Sigma^0,\Z)$.
It is given by 
$$\bigoplus_{i=0}^r (1_{U(1)}\times\ori_{\partial \bar U_i})\Z\oplus \bigoplus_{i=0}^r (\ori_{U(1)}\times 1_{\partial \bar U_i})\Z\rightarrow \Z\oplus \bigoplus_{i=0}^r(\ori_{U(1)}\times 1_{\partial \bar U_i})\Z\ ,$$
where the first component
maps
$\sum_{i=0}^r a_i  (1_{U(1)}\times\ori_{\partial \bar U_i})$ to $\sum_{i=0}^r
a_i$, and the second component is the identity.
Let
$$\alpha_1: (\ori_{U(1)}\times 1_{\Sigma^0})\Z\oplus   \bigoplus_{i=0}^rH^1([U(1)/_{\chi_i}\Gamma_i],\Z)\rightarrow \Z\oplus \bigoplus_{i=0}^r(\ori_{U(1)}\times 1_{\partial \bar U_i})\Z
$$
be the induced map.
We have
$$\alpha_1(\ori_{U(1)}\times 1_{\Sigma^0})=0\oplus \oplus_{i=0}^r (\ori_{U(1)}\times 1_{\partial \bar U_i})\ .$$
We now describe the restriction of $\alpha_1$ to the summand
$H^1([U(1)/_{\chi_i}\Gamma_i],\Z)$.
It is given by the composition of pull-backs along the following sequence of maps:
\begin{eqnarray*}U(1)\times \partial \bar U_i\cong [U(1)\times \partial \tilde{\bar U}_i/_1\Gamma_i]\stackrel{I_{\phi_i}}{\rightarrow} [U(1)\times \partial \tilde{\bar U}_i/_{\chi_i}\Gamma_i]\rightarrow
[U(1)\times \tilde{\bar U}_i/_{\chi_i}\Gamma_i]&&\\\rightarrow [U(1)/_{\chi_i}\Gamma_i]&&\ ,\end{eqnarray*}
where $I_{\phi_i}$ is induced by the map $\phi_i$ (see \ref{con})
$I_{\phi_i}(z,\tilde s):=(\phi_i(\tilde s)z,\tilde s)$, and the remaining maps are the obvious inclusions and projections. In the case $i=0$ we set $\phi_0:=u$. 

We have
$H^1([U(1)/_{\chi_i}\Gamma_i],\Z)\cong H^1(U(1)\times_{\Gamma_i,\chi_i} E\Gamma_i,\Z)$. We consider the $U(1)$-bundle
$U(1)\times_{\Gamma_i,\chi_i} E\Gamma_i\rightarrow B\Gamma_i$.
Using the Serre spectral sequence we see that the restriction to the fibre
$r^*$ fits into an 
exact sequence
$$0\rightarrow H^1([U(1)/_{\chi_i}\Gamma_i],\Z)\stackrel{r^*}{\rightarrow} \Z
\stackrel{\chi_i}{\rightarrow}\Z/n_i\Z\ .$$
Similarly, restriction to the fibre of the bundle
$(U(1)\times \partial \tilde{\bar U}_i)_{\Gamma_i,1}\times E\Gamma_i\rightarrow B\Gamma_i$ gives an exact  sequence
$$0\rightarrow H^1([U(1)\times \partial \tilde{\bar U}_i/_1\Gamma_i],\Z)\rightarrow \Z\oplus \Z \stackrel{\pr_2}{\rightarrow} \Z/n_i\Z\ ,$$
where we use the basis
$\Z\oplus \Z\cong (\ori_{U(1)}\times 1_{\partial \tilde{\bar U}_i})\Z\oplus(
1_{U(1)}\times \ori_{\partial \tilde{\bar U}_i})\Z$.

Let $a\in \Z$ represent an element of $H^1([U(1)/_{\chi_i}\Gamma_i],\Z)$,
i.e. $\chi_i [a]=0\in \Z/n_i\Z$.
Then one can check that $\alpha_1(a)=(a,\deg(\phi_i)a)$.
Fortunately, as observed in \ref{con}, $[\deg(\phi_i)]=\chi$
in $\hat \Gamma_i\cong \Z/n_i\Z$ so that $n_i|   \deg(\phi_i)a$, and thus
$( a,\deg(\phi_i)a)\in H^1([U(1)\times \partial \tilde{\bar U}_i/_1\Gamma_i],\Z)$.
Combining these calculations we obtain the following explicit description of $$\alpha_1:\Z\oplus \bigoplus_{i=0}^r \ker(\chi_i:\Z\rightarrow \Z/n_i\Z)\rightarrow \Z\oplus \Z^{r+1} \ ,$$
$$\alpha_1(x,(a_0,\dots,a_r))=(\sum_{i=0}^r \frac{\deg(\phi_i)a_i}{n_i},(a_0+x,\dots,a_r+x))\ ,$$
where on the right-hand side we identify
$\Z^{r+1}\cong \bigoplus_{i=0}^r (\ori_{U(1)}\times 1_{\partial \bar U_i})\Z$.

We now have collected sufficient information on $H^2(E,\Z)$ in order to calculate $c_1(E)$. By the compatibility of the Mayer-Vietoris sequences with the pull-back $\pi^*:H^2(B,\Z)\rightarrow H^2(E,\Z)$ we get the diagram
$$\begin{array}{ccccccccc}
0&\rightarrow&\coker(\alpha_1)&\rightarrow&H^2(E,\Z)&\rightarrow&\Z^{2g}\oplus\bigoplus_{i=1}^r \hat \Gamma_i/\chi_i\hat \Gamma_i&\rightarrow&0\\
&&f\uparrow &&\uparrow&&g\uparrow&&\\
0&\rightarrow&\Z&\rightarrow&H^2(B,\Z)&\stackrel{\oplus_{i=0}^rt_i^*}{\rightarrow}&\bigoplus_{i=0}^r \hat \Gamma_i&\rightarrow&0\end{array}\ ,$$
where $t_i:[p_i/\Gamma_i]\rightarrow B$ is the canonical embedding.
We must determine generators of $\ker(f)$ and $\ker(g)$.
We have a factorization of $f$ as
$\Z\stackrel{(\id,0)}{\rightarrow} \Z\oplus \Z^{r+1}\rightarrow \coker(\alpha_1)$. We see that
$f(b)=0$ is equivalent to the condition that the system
\begin{eqnarray*}
b&=&x \sum_{i=0}^r \frac{\deg(\phi_i)}{n_i}\\
\chi_i [x]&=&0 \in \Z/n_i\Z\ ,\quad i=0,\dots r
\end{eqnarray*}
has a solution $x\in \Z$. We see that
$\ker(f)\subset \Z$ is a non-trivial subgroup, and we fix the generator
$e\in \Z$ which is given by the component of $c_1(E)$. It is determined by the subgroup up to sign.

The kernel of $g$ is the sum of the kernels of the projections
$\hat \Gamma_i\rightarrow \hat \Gamma_i/\chi_i\hat \Gamma_i$.
In order to find the generators which correspond to the Chern character of $E$
we use the fact that the Chern character is compatible with restriction.
We consider the pull-back
$$\begin{array}{ccc}
E_{p_i}&\rightarrow &E\\
\downarrow&&\downarrow\\
{}[p_i/\Gamma_i]&\stackrel{t_i}{\rightarrow}&B\end{array} .$$
Since we know that
$c_1(E_{p_i})=\chi_i\in \hat \Gamma_i$ we see that
$t_i^* c_1(E)=\chi_i$ is the correct generator of the kernel of the corresponding component of $g$.

Combining these calculations we get
$$c_1(E)=(e,(\chi_1,\dots,\chi_1))\in \Z\oplus \bigoplus_{i=1}^r \hat \Gamma_i\ ,$$
where $e$ was described above.

\subsubsection{}

We now compute
$H^3(E,\Z)$, again using a Mayer-Vietoris sequence.
Let $[U(1)/_{\chi_i} \Gamma_i]$ denote the orbispace given by the action of $\Gamma_i$ on $U(1)$ via $\chi_i$.
The relevant part of the Mayer-Vietoris sequence has the form
\begin{eqnarray*}
 H^2(U(1)\times \Sigma^0,\Z)\oplus\bigoplus_{i=0}^r
H^2([U(1)/_{\chi_i}\Gamma_i],\Z)\rightarrow \bigoplus_{i=0}^r H^{2}(U(1)\times \partial \bar U_i,\Z)&&\\\rightarrow
H^3(E,\Z)\rightarrow \bigoplus_{i=0}^r
H^3([U(1)/_{\chi_i}\Gamma_i],\Z)\rightarrow 0\
.\end{eqnarray*}
We now use the facts that the restriction
$H^2([U(1)/_{\chi_i}\Gamma_i],\Z)\rightarrow H^{2}(U(1)\times \partial \bar U_i,\Z)$ is trivial, that the cokernel of $H^2(U(1)\times \Sigma^0,\Z)\rightarrow
\bigoplus_{i=0}^r H^{2}(U(1)\times \partial \bar U_i,\Z)$ is isomorphic to $\Z$, and that $H^3([U(1)/_{\chi_i}\Gamma_i],\Z)\cong \Ann(\chi_i)$,
where the definition of $\Ann(\chi_i)\subset \hat \Gamma_i$ uses the
ring structure of $\hat  \Gamma_i$ (see \ref{zwde} for the computation
of $H^3([U(1)/_{\chi_i}\Gamma_i],\Z)$).
The sequence thus simplifies to
$$0\rightarrow \Z\stackrel{\delta_E}{\rightarrow} H^3(E,\Z)\rightarrow
\bigoplus_{i=1}^r \Ann(\chi_i)\rightarrow 0\ .$$
Let $\pi:E\rightarrow B$ be the projection. Then the
following diagram commutes:
$$\begin{array}{ccc}
\Z&\stackrel{\delta_E}{\rightarrow}&H^3(E,\Z)\\
\|&&\pi_!\downarrow\\
\Z&\stackrel{\delta}{\rightarrow}&H^2(B,\Z)
\end{array}
\ .$$
Therefore the decomposition  $H^2(B,\Z)=\Z\oplus \bigoplus_{i=1}^r \hat
\Gamma_i$ induces a split $s_E:H^3(E,\Z)\rightarrow \Z$,
so that we obtain an identification
$$H^3(E,\Z)\cong \Z\oplus \bigoplus_{i=1}^r
\Ann(\chi_i)\ .$$ Note that this decomposition is again canonical.
A cohomology class $h\in H^3(E,\Z)$ is thus identified with an element
$$(f,(a_1,\dots,a_r))\in  \Z\oplus \Ann(\chi_1)\oplus
\dots\oplus \Ann(\chi_r)\ .$$

\subsubsection{}

It follows from Proposition \ref{ttz} that
the topological type of $E$ is classified by $c_1(E)$.

We further observe that
$\pi_!:H^3(E,\Z)\rightarrow H^2(B,\Z)$ is injective. Therefore we can characterize a class in $H^3(E,\Z)$ by its image under $\pi_!$. It follows that automorphisms of  the $U(1)$-bundle $E$ act trivially on $H^3(E,\Z)$.
We see that the  isomorphism class of the pair $(E,h)$ is determined by
$$(c_1(E),\pi_!(h))=(e,(\chi_1,\dots,\chi_r),f,(a_1,\dots,a_r))\in H^2(B,\Z)\oplus H^2(B,\Z)$$
(see \ref{h2b} for the notation).
It therefore makes sense to calculate the $T$-dual pair $(\hat E,\hat h)$ in terms of its topological invariants $(c_1(\hat E),\hat \pi_!(\hat h))$.
We get
$$(c_1(\hat E),\hat \pi_!(\hat h))=(-f,(-a_1,\dots,-a_r),-e,(-\chi_1,\dots,-\chi_r))\ .$$

\end{document}